\documentclass[11pt]{amsart}
\title{Derived Recollements and Generalised AR Formulas} 
\author{Samuel Dean and Jeremy Russell}
\address{University of Manchester, Manchester, UK}
\email{samuel.dean-3@postgrad.manchester.ac.uk}
\address{The College of New Jersey, Ewing, New Jersey}
\email{russelj1@tcnj.edu}
\usepackage{amsmath,amssymb,amsthm, amsfonts, xypic, graphicx}

\usepackage[hidelinks=true]{hyperref}
\usepackage[all]{xy}

\theoremstyle{definition}
\newtheorem{thm}{Theorem}
\newtheorem{prop}[thm]{Proposition}
\newtheorem{lem}[thm]{Lemma}
\newtheorem{cor}[thm]{Corollary}

\newenvironment{pf}{\paragraph{{\sc \dbf{Proof}}}}{\par\vspace{1cm}}
\theoremstyle{definition}
\newtheorem{defn}{\underline{Definition}}
\newtheorem*{rem}{Remark}
\newtheorem*{thm*}{Theorem}
\newtheorem*{prop*}{Proposition}

\newtheorem*{com}{Comment}
\newtheorem*{lem*}{Lemma}
\newtheorem*{cor*}{Corollary}
\newtheorem*{diag*}{AGJ Adjunction Diagram}
\newtheorem*{agjr*}{The Auslander-Gruson-Jensen Recollement}
\newtheorem*{defr*}{The Defect Recollement}
\newtheorem*{evr*}{The Evaluation Recollement}
\newtheorem*{comp*}{Right Derived Dual Complex}

\newtheorem{diag}{Diagram}

\newtheorem*{ex*}{Example}

\renewcommand{\qed}{\blacksquare}

\newcommand{\newdual}{\textsf{T}}

\newcommand{\la}{\lambda}
\newcommand{\ra}{\rho}
\newcommand{\lla}{\mu}

\newcommand{\set}[1]{\left\{#1\right\}}
\newcommand{\blank}{\hspace{0.05cm}\underline{\ \ }\hspace{0.1cm} }
\newcommand{\ev}{\textsf{ev}}
\renewcommand{\mod}{\textsf{mod}}
\newcommand{\Mod}{\textsf{Mod}}
\newcommand{\A}{\mathcal{A}}
\newcommand{\B}{\mathcal{B}}
\newcommand{\C}{\mathcal{C}}
\newcommand{\D}{\mathcal{D}}
\newcommand{\E}{\mathcal{E}}

\renewcommand{\P}{\mathcal{P}}

\newcommand{\Y}{\textsf{Y}}

\newcommand{\sub}{\subseteq}
\newcommand{\tfp}{\textsf{pp}}
\newcommand{\fp}{\textsf{fp}}

\newcommand{\dbf}[1]{\textbf{\text{#1}}}

\usepackage{pgfpages}
\usepackage{tikz}
\usetikzlibrary{matrix,arrows,snakes,backgrounds}

\renewcommand{\:}{\colon}

\newcommand{\eps}{\varepsilon}
\newcommand{\To}{\longrightarrow}

\newcommand{\ab}{\textsf{Ab}}

\newcommand{\im}{\textsf{Im}}

\newcommand{\coker}{\textsf{Coker}}
\renewcommand{\ker}{\textsf{Ker}}

\newcommand{\tor}{\textsf{Tor}}
\newcommand{\ext}{\textsf{Ext}}

\newcommand{\tr}{\textrm{Tr}}

\newcommand{\nat}{\textsf{Nat}}
\renewcommand{\hom}{\textsf{Hom}}

\newcommand{\res}{\textsf{R}}
\newcommand{\h}{\textsf{h}}
\newcommand{\agj}{\textsf{D}_\textsf{A}}
\newcommand{\agr}{\textsf{D}_\textsf{R}}
\newcommand{\agl}{\textsf{D}_\textsf{L}}
\newcommand{\agd}{\textsf{D}}
\newcommand{\extd}{\textsf{W}}
\begin{document}
\maketitle

\begin{abstract} The Defect Recollement, Restriction Recollement, Auslander-Gruson-Jensen Recollement, and others, are shown to be instances of a general construction using derived functors and methods from stable module theory.  The right derived functors $\textsf{W}_k:=R_k(\hspace{0.05cm}\underline{\ \ }\hspace{0.1cm} )^*$ are computed and it is shown that the functor $\textsf{W}_2:=R_2(\hspace{0.05cm}\underline{\ \ }\hspace{0.1cm} )^*$ is right exact and restricts to a duality $\textsf{W}$ of the defect zero functors.   The duality $\textsf{W}$ satisfies two identities which we call the Generalised Auslander-Reiten formulas. We show that $\textsf{W}$ restricts to the generalised Auslander-Bridger transpose and show that the Generalised Auslander-Reiten formulas reduce to the well-known Auslander-Reiten formulas. 

 \end{abstract}

\setcounter{tocdepth}{1}
\tableofcontents

\section{Introduction}
In \cite{dr1}, the authors studied the Auslander-Gruson-Jenson functor $$\agj\:\fp(\Mod(R),\ab)\To (\mod(R^{op}),\ab)$$ which is an exact contravariant functor sending representable functors $(M,\blank)$ to tensor functors $M\otimes\blank$, showing that the Auslander-Gruson-Jensen functor admits both a left and right adjoint, thus establishing that $\agj$ is part of a recollement called the Auslander-Gruson-Jensen recollement. In this paper, we streamline and generalise the proof of these assertions.

In \cite{kraus}, Krause shows that for a finitely accessible category $\A$, the restriction functor $$\res\:\fp(\A^{op},\ab)\To ((\fp \A)^{op},\ab)$$ admits two adjoints.   In the case that $\A=\Mod(R)$, then the description of the both kernels of $\agj$ and $\res$ is the full subcategory of functors arising from pure exact sequences in $\Mod(R)$.  The goal of this paper is to fully explain the connection between the functors $\agj$ and $\res$, both by realizing them as pieces of a more general construction, and giving a relationship between them which explains why they have dual kernels.   The main contributions of the paper are the following.

%A reflective adjunction is an adjoint pair $\lambda\dashv \ra$ where $\ra$ is fully faithful.   Under reasonable conditions we can construct recollements from reflective adjunctions by using 0-th derived functors.   

\begin{thm*}[Derived Recollement Theorem]\label{recthm}Suppose that $\C,\D$ are both abelian categories with enough projectives and we have an adjunction

\begin{center}\begin{tikzpicture}
\matrix (m) [ampersand replacement= \&,matrix of math nodes, row sep=3em, 
column sep=4em,text height=1.5ex,text depth=0.25ex] 
{\C\&\D\\}; 
\path[->,thick, font=\scriptsize]
(m-1-1) edge[bend left=50] node[auto]{$\la$}(m-1-2)
(m-1-2) edge[bend left=50] node[below]{$\ra$}(m-1-1);
\end{tikzpicture}\end{center} Assume that 
\begin{enumerate}
\item For any projective $P\in \C$ there is some object $D\in \D$ such that $P\cong\ra (D)$.
\item $\rho$ sends projectives to projectives.
\end{enumerate}  If  the adjunction $\la\dashv\ra$ is also reflective, meaning that $\ra$ is fully faithful, then it can be extended to a recollement

\begin{center}\begin{tikzpicture}
\matrix (m) [ampersand replacement= \&,matrix of math nodes, row sep=5em, 
column sep=5em,text height=1.5ex,text depth=0.25ex] 
{\ker(\la)\&\C\&\D\\}; 
\path[->,thick, font=\scriptsize]
(m-1-1) edge node[above]{$i$}(m-1-2)
(m-1-2) edge[bend right=50] node[above]{$i_L$}(m-1-1)
(m-1-2) edge[bend left=50] node[below]{$i_R$}(m-1-1)
(m-1-3)edge[bend left=50] node[below]{$\ra$}(m-1-2)
(m-1-3) edge[bend right=50] node[above]{$L_0\ra$}(m-1-2)
(m-1-2) edge node[above]{$\la$}(m-1-3);
\end{tikzpicture}\end{center} where $L_0\ra$ is the zeroth left derived functor.
In particular, $\la$ is exact.
\end{thm*}

We apply this result to produce various examples of recollements.  One of our examples occurs if $\C$ is a finite $k$-variety with weak kernels and weak cokernels such that $\mod(\C)$ is a dualising $k$-variety.  In this case we show that $\C$ is also a dualising $k$-variety. This generalises \cite[Proposition 2.5]{dualvar}.   We also show that our motivating examples, namely the Auslander-Gruson-Jensen functor $\agj$ and the Restriction functor $\res$,  are both localization functors appearing in derived recollements.

Our focus then shifts to fully understanding the connection between these two recollements.  This is ultimately achieved by looking at the functors \begin{center}\begin{tikzpicture}
\matrix (m) [ampersand replacement= \&,matrix of math nodes, row sep=3em, 
column sep=5em,text height=1.5ex,text depth=0.25ex] 
{\fp(\A,\ab)\&\fp(\A^{op},\ab)\\}; 
\path[->,thick, font=\scriptsize]
(m-1-1) edge[bend left=50] node[auto]{$(\blank)^*$}(m-1-2)
(m-1-2) edge[bend left=50] node[below]{$(\blank)^*$}(m-1-1);
\end{tikzpicture}\end{center} which are completely determined by contravariant left exactness and the formulas $$(X,\blank)^*\cong (\blank, X)$$ $$(\blank,X)^*\cong (X,\blank)$$  In \cite{defect} the left derived functors $L^k(\blank)^*$ are computed and it is determined that the only non-trivial left derived functor is an exact functor $L^0(\blank)^*$.  In this paper, we complete the description of the derived functors of $(\blank)^*$ by computing the right derived functors $\extd_k:=R_k(\blank)^*$.  In doing so we recover the following.

%This naturally leads us to the categories $\fp_0(\Mod(R),\ab)$ and $\fp_0(\Mod(R)^{op},\ab)$ consisting of all functors in the categories $\fp(\Mod(R),\ab)$ and $\fp_0(\Mod(R)^{op},\ab)$ arising from pure exact sequences.  

\begin{thm*}\label{defectzerodual}For any abelian category $\A$, the functors $\extd_2=R_2(\blank)^*$ form a pair of right exact functors \begin{center}\begin{tikzpicture}
\matrix (m) [ampersand replacement= \&,matrix of math nodes, row sep=5em, 
column sep=4em,text height=1.5ex,text depth=0.25ex] 
{\fp(\A,\ab)\&\fp(\A^{op},\ab)\\}; 
\path[->,thick, font=\scriptsize]
(m-1-1)edge[bend right=50] node[below]{$\extd_2$}(m-1-2)
(m-1-2) edge[bend right=50] node[above]{$\extd_2$}(m-1-1);
%(m-1-2) edge[bend left] node[below]{$\agr$}(m-1-1);
\end{tikzpicture}\end{center}  that restricts to a duality $\extd$

\begin{center}\begin{tikzpicture}
\matrix (m) [ampersand replacement= \&,matrix of math nodes, row sep=5em, 
column sep=4em,text height=1.5ex,text depth=0.25ex] 
{\fp_0(\A,\ab)\&\fp_0(\A^{op},\ab)\\}; 
\path[->,thick, font=\scriptsize]
(m-1-1)edge[bend right=50] node[below]{$\extd$}(m-1-2)
(m-1-2) edge[bend right=50] node[above]{$\extd$}(m-1-1);
\end{tikzpicture}\end{center}

\end{thm*}

This theorem is then applied to recover two formulas that look strikingly similar to the Auslander-Reiten formulas. At the end of the paper, we prove that they do in fact reduce to the original Auslander-Reiten formulas.

\begin{thm*}[The Generalised Auslander-Reiten Formulas]Suppose that $\A$ is abelian category with enough projectives and injectives.  Then we have
\begin{eqnarray}
\underline{\hom}(\blank, C)&\cong& \extd \ext^1(C,\blank)\\
\extd \ext^1(\blank, C)&\cong&\overline{\hom}(C,\blank)
\end{eqnarray}
\end{thm*}

These two formulas allow us to achieve our ultimate goal of understanding how the Auslander-Gruson-Jensen Recollement and Restriction Recollement relate to each other.  The duality $\extd$ together with the localization properties of $\agj$ and $\res$ are used to show the following.

\begin{thm*}The duality $\extd$ induces a duality $\newdual$ making the following diagram of functors commute.

\begin{center}\begin{tikzpicture}
\matrix(m)[ampersand replacement=\&, matrix of math nodes, row sep=2em, column sep=2em, text height=1.5ex, text depth=0.25ex]
{\fp_0(\Mod(R),\ab)\&\fp_0(\Mod(R)^{op},\ab)\&\fp_0(\Mod(R),\ab)\\
(\underline{\mod}(R^{op}),\ab)\&(\underline{\mod}(R)^{op},\ab)\&(\underline{\mod}(R^{op}),\ab)\\}; 
\path[->,thick, font=\scriptsize]
(m-1-1) edge node[above]{$\extd$}(m-1-2)
(m-1-1) edge node[left]{$\agj$}(m-2-1)
(m-1-2) edge node[above]{$\extd$}(m-1-3)
(m-1-2) edge node[right]{$\res$}(m-2-2)
(m-1-3) edge node[right]{$\agj$}(m-2-3)
(m-2-1) edge node[below]{$\newdual$}(m-2-2)
(m-2-2) edge node[below]{$\newdual$}(m-2-3);
\end{tikzpicture}\end{center} The duality $\newdual$ is isomorphic to the generalized Auslander-Bridger transpose $\tr_*$ studied by Herzog in \cite{herzogcontra}.  Hence the duality $\tr_*$ can be recovered via a universal property.
\end{thm*}

The paper is organized as follows.  In Section \ref{section2}, we recall the functor categories $\fp(\A,\ab)$ and $(\mod(R^{op}),\ab)$ and describe some of their well known properties.   The category $\fp(\A,\ab)$ is the category of the finitely presented functors $F\:\A\to \ab$ together with the natural transformations between them.  A functor $F$ is finitely presented if there exists an exact sequence  \begin{center}\begin{tikzpicture}
\matrix(m)[ampersand replacement=\&, matrix of math nodes, row sep=3em, column sep=2em, text height=1.5ex, text depth=0.25ex]
{(Y,\blank)\&(X,\blank)\&F\&0\\}; 
\path[->,thick, font=\scriptsize]
(m-1-1) edge node[auto]{}(m-1-2)
(m-1-2) edge node[auto]{}(m-1-3)
(m-1-3) edge node[auto]{}(m-1-4);
\end{tikzpicture}\end{center}    

The Yoneda embedding $\Y\:\A\To \fp(\A,\ab)$ is the fully faithful contravariant functor defined by $\Y(X)=(X,\blank)$.   In \cite{coh}, Auslander constructs an exact contravariant functor $w\:\fp(\A,\ab)\To \A$ satisfying the property that $w(X,\blank)\cong X$.  He further shows that the defect is an adjoint to the Yoneda embedding.  For each finitely presented functor he constructs a four-term exact sequence called the \dbf{defect sequence}: 
\begin{center}\begin{tikzpicture}
\matrix(m)[ampersand replacement=\&, matrix of math nodes, row sep=3em, column sep=2.5em, text height=1.5ex, text depth=0.25ex]
{0\&F_0\&F\&\big(w(F),\blank\big)\&F_1\&0\\}; 
\path[->,thick, font=\scriptsize]
(m-1-1) edge node[auto]{}(m-1-2)
(m-1-2) edge node[auto]{}(m-1-3)
(m-1-3) edge node[auto]{$\varphi$}(m-1-4)
(m-1-4) edge node[auto]{}(m-1-5)
(m-1-5) edge node[auto]{}(m-1-6);
\end{tikzpicture}\end{center}

The category $(\mod(R^{op}),\ab)$ consists of all functors $F\:\mod(R^{op})\to \ab$ together with the natural transformations between them.   The tensor embedding is the fully faithful right exact functor $t\:\Mod(R)\to (\mod(R^{op}),\ab)$ given by $t(M)=M\otimes\blank$.  In \cite{gj}, Gruson and Jensen show that evaluation at the ring $$\ev_R\:(\mod(R^{op}),\ab)\To \Mod(R)$$ is the right adjoint to $t$.  

The Auslander-Gruson-Jensen functor is the functor $$\agj\:\fp(\Mod(R),\ab)\to (\mod(R^{op}),\ab)$$ defined by $$\agj:=R_0(t\circ w)$$ An alternate description of $\agj$ is given by its exactness and the fact that it sends any representable $(X,\blank)$ to the tensor functor $X\otimes\blank$.  The functor $\agj$ was studied by the authors in \cite{dr1} where it is shown to have both adjoints.

Finally, the definition of the the stable categories $\underline{\A}$ and $\overline{\A}$ are reviewed at the end of the section.   The category $\underline{\A}$ is the quotient category obtained by factoring out the ideal of all morphisms in $\A$ which factor through a projective. These categories will be used later on in the paper and so we have briefly included their definitions.

In Section \ref{section3}, we introduce recollements. A recollement is a diagram of functors between abelian categories $\A,\B,\C$,

\begin{center}\begin{tikzpicture}
\matrix (m) [ampersand replacement= \&,matrix of math nodes, row sep=5em, 
column sep=4em,text height=1.5ex,text depth=0.25ex] 
{\A\&\B\&\C\\}; 
\path[->,thick, font=\scriptsize]
(m-1-1) edge node[above]{$i$}(m-1-2)
(m-1-2) edge[bend right=50] node[above]{$i_L$}(m-1-1)
(m-1-2) edge[bend left=50] node[below]{$i_R$}(m-1-1)
(m-1-3)edge[bend left=50] node[below]{$f_R$}(m-1-2)
(m-1-3) edge[bend right=50] node[above]{$f_L$}(m-1-2)
(m-1-2) edge node[above]{$f$}(m-1-3);
\end{tikzpicture}\end{center}
where $i_L\dashv i\dashv i_R$, $f_L\dashv f\dashv f_R$, the functors $i$, $f_L$ and $f_R$ are fully faithful, and the essential image of $i$ is the kernel of $f$. We generalise Auslander's method for constructing the defect sequence in \cite{coh}. We introduce reflective and coreflective adjunctions. A reflective (respectively, coreflective) adjunction is an adjoint pair $\la\dashv\ra$ where $\ra$ (respectively, $\la$) is fully faithful, and explain that any localisation of abelian categories (i.e. an exact functor $\lambda$ which is part of a reflective adjunction $\la\dashv\ra$) is a Serre localisation. We then prove one of our main results: That any reflective adjunction $\la\dashv\ra$ where $\ra$ which maps projectives to projectives and maps at least one object onto every projective is part of a recollement $\mu=L_0\ra\dashv\la\dashv\ra$. We call this result the Derived Recollement Theorem.

In Section \ref{section4}, we introduce the notion of a derived recollement. This is simply a recollement gotten by an application of the Derived Recollement Theorem. We then give several examples: 
\begin{enumerate}
\item If $\A$ is abelian and has enough projectives then there is a recollement \begin{center}\begin{tikzpicture}
\matrix (m) [ampersand replacement= \&,matrix of math nodes, row sep=5em, 
column sep=6em,text height=1.5ex,text depth=0.25ex] 
{\ker(w)\&\fp(\A^{op},\ab)\&\A\\}; 
\path[->,thick, font=\scriptsize]
(m-1-1) edge node[above]{$i$}(m-1-2)
(m-1-2) edge[bend right=50] node[above]{$(\blank)^0$}(m-1-1)
(m-1-2) edge[bend left=50] node[below]{$(\blank)_0$}(m-1-1)
(m-1-3)edge[bend left=50] node[below]{$\Y$}(m-1-2)
(m-1-3) edge[bend right=50] node[above]{$L_0\Y$}(m-1-2)
(m-1-2) edge node[above]{$w$}(m-1-3);
\end{tikzpicture}\end{center} 

\item If $\C$ is finitely accessible and cocomplete and $v:((\fp \C)^{op},\ab)\to C$ is the unique colimit preserving functor which sends a functor $(-,C)$ to the corresponding object $C\in\C$. There is a recollement \begin{center}\begin{tikzpicture}
\matrix (m) [ampersand replacement= \&,matrix of math nodes, row sep=5em, 
column sep=6em,text height=1.5ex,text depth=0.25ex] 
{\ker(v)\&((\fp\C)^{op},\ab)\&\C\\}; 
\path[->,thick, font=\scriptsize]
(m-1-1) edge node[above]{$i$}(m-1-2)
(m-1-2) edge[bend right=50] node[above]{$(\blank)^0$}(m-1-1)
(m-1-2) edge[bend left=50] node[below]{$(\blank)_0$}(m-1-1)
(m-1-3)edge[bend left=50] node[below]{$\Y$}(m-1-2)
(m-1-3) edge[bend right=50] node[above]{$L_0\Y$}(m-1-2)
(m-1-2) edge node[above]{$v$}(m-1-3);
\end{tikzpicture}\end{center} 
\item Let $R$ be any ring. There is a recollement \begin{center}\begin{tikzpicture}
\matrix (m) [ampersand replacement= \&,matrix of math nodes, row sep=5em, 
column sep=6em,text height=1.5ex,text depth=0.25ex] 
{\ker(\ev_R)\&(\mod(R^{op}),\ab)\&\Mod(R)\\}; 
\path[->,thick, font=\scriptsize]
(m-1-1) edge node[above]{$i$}(m-1-2)
(m-1-2) edge[bend right=50] node[above]{$(\blank)^0$}(m-1-1)
(m-1-2) edge[bend left=50] node[below]{$(\blank)_0$}(m-1-1)
(m-1-3)edge[bend left=50] node[below]{$R_0t$}(m-1-2)
(m-1-3) edge[bend right=50] node[above]{$t$}(m-1-2)
(m-1-2) edge node[above]{$\ev_R$}(m-1-3);
\end{tikzpicture}\end{center}
\item Let $\C$ be an additive category with weak cokernels. There is a recollement \begin{center}\begin{tikzpicture}
\matrix (m) [ampersand replacement= \&,matrix of math nodes, row sep=5em, 
column sep=6em,text height=1.5ex,text depth=0.25ex] 
{\ker(\ev_\C)\&\fp(\mod(\C^{op}),\ab)\&\mod(\C)\\}; 
\path[->,thick, font=\scriptsize]
(m-1-1) edge node[above]{$i$}(m-1-2)
(m-1-2) edge[bend right=50] node[above]{$(\blank)^0$}(m-1-1)
(m-1-2) edge[bend left=50] node[below]{$(\blank)_0$}(m-1-1)
(m-1-3)edge[bend left=50] node[below]{$R_0t$}(m-1-2)
(m-1-3) edge[bend right=50] node[above]{$t$}(m-1-2)
(m-1-2) edge node[above]{$\ev_\C$}(m-1-3);
\end{tikzpicture}\end{center} 
\item The Auslander-Gruson-Jensen Recollement. For any ring $R$, there is a recollement
\begin{center}\begin{tikzpicture}
\matrix (m) [ampersand replacement= \&,matrix of math nodes, row sep=5em, 
column sep=3em,text height=1.5ex,text depth=0.25ex] 
{\ker(\agj)\&\fp(\Mod(R),\ab)\&(\mod(R^{op}),\ab)^{op}\\}; 
\path[->,thick, font=\scriptsize]
(m-1-1) edge node[above]{$i$}(m-1-2)
(m-1-2) edge[bend right=50] node[above]{$(\blank)^p$}(m-1-1)
(m-1-2) edge[bend left=50] node[below]{$(\blank)_q$}(m-1-1)
(m-1-3)edge[bend left=50] node[below]{$\agr$}(m-1-2)
(m-1-3) edge[bend right=50] node[above]{$\agl=L_0(\agr)$}(m-1-2)
(m-1-2) edge node[above]{$\agj$}(m-1-3);
\end{tikzpicture}\end{center} 
\item The Restriction Recollement. For any finitely accessible category $\C$, there is a recollement
\begin{center}\begin{tikzpicture}
\matrix (m) [ampersand replacement= \&,matrix of math nodes, row sep=5em, 
column sep=3em,text height=1.5ex,text depth=0.25ex] 
{\ker(\res)\&\fp(\C^{op},\ab)\&((\fp\C)^{op},\ab)^{op}\\}; 
\path[->,thick, font=\scriptsize]
(m-1-1) edge node[above]{$i$}(m-1-2)
(m-1-2) edge[bend right=50] node[above]{$(\blank)^0$}(m-1-1)
(m-1-2) edge[bend left=50] node[below]{$(\blank)_0$}(m-1-1)
(m-1-3)edge[bend left=50] node[below]{$\overleftarrow{-}$}(m-1-2)
(m-1-3) edge[bend right=50] node[above]{$L_0(\overleftarrow{-})$}(m-1-2)
(m-1-2) edge node[above]{$\res$}(m-1-3);
\end{tikzpicture}\end{center}  
\end{enumerate}

Let $\A$ be an abelian category. In Section \ref{section5}, we construct the duality $\extd$ by computing the right derived functors of the dual $(\blank)^*$  \begin{center}\begin{tikzpicture}
\matrix (m) [ampersand replacement= \&,matrix of math nodes, row sep=3em, 
column sep=5em,text height=1.5ex,text depth=0.25ex] 
{\fp(\A,\ab)\&\fp(\A^{op},\ab)\\}; 
\path[->,thick, font=\scriptsize]
(m-1-1) edge[bend left=50] node[auto]{$(\blank)^*$}(m-1-2)
(m-1-2) edge[bend left=50] node[below]{$(\blank)^*$}(m-1-1);
\end{tikzpicture}\end{center} which were first defined by Fisher-Palmquist and Newell in their paper \cite{fisheradjoint} on adjoint pairs and bifunctors.

Let $F$ be any finitely presented functor and take any projective resolution.

\begin{center}\begin{tikzpicture}
\matrix(m)[ampersand replacement=\&, matrix of math nodes, row sep=3em, column sep=2.5em, text height=1.5ex, text depth=0.25ex]
{0\&(Z,\blank)\&(Y,\blank)\&(X,\blank)\&F\&0\\}; 
\path[->,thick, font=\scriptsize]
(m-1-1) edge node[auto]{}(m-1-2)
(m-1-2) edge node[auto]{$(g,\blank)$}(m-1-3)
(m-1-3) edge node[auto]{$(f,\blank)$}(m-1-4)
(m-1-4) edge node[auto]{}(m-1-5)
(m-1-5) edge node[auto]{}(m-1-6);
\end{tikzpicture}\end{center}  The functors $\extd_k(F)$ coincide with the homology of the complex

\begin{center}\begin{tikzpicture}
\matrix(m)[ampersand replacement=\&, matrix of math nodes, row sep=3em, column sep=2.5em, text height=1.5ex, text depth=0.25ex]
{0\&(\blank,X)\&(\blank,Y)\&(\blank,Z)\&0\\}; 
\path[->,thick, font=\scriptsize]
(m-1-1) edge node[auto]{}(m-1-2)
(m-1-2) edge node[auto]{$(\blank,f)$}(m-1-3)
(m-1-3) edge node[auto]{$(\blank,g)$}(m-1-4)
(m-1-4) edge node[auto]{}(m-1-5);
\end{tikzpicture}\end{center}  In particular, we show that $\extd_0(F)\cong(F)^*$ and $\extd_2(F)\cong\coker(\blank, g)$.  Hence given the presentation of $F$ above we get the following presentation of $\extd_2(F)$.

\begin{center}\begin{tikzpicture}
\matrix(m)[ampersand replacement=\&, matrix of math nodes, row sep=3em, column sep=2.5em, text height=1.5ex, text depth=0.25ex]
{0\&(\blank, V)\&(\blank,Y)\&(\blank,Z)\&\extd_2(F)\&0\\}; 
\path[->,thick, font=\scriptsize]
(m-1-1) edge (m-1-2)
(m-1-2) edge node[above]{$(\blank, m)$}(m-1-3)
(m-1-3) edge node[auto]{$(\blank,g)$}(m-1-4)
(m-1-4) edge node[auto]{}(m-1-5)
(m-1-5) edge node[auto]{}(m-1-6);
\end{tikzpicture}\end{center}

This is used to show that $\extd_2\circ\extd_2(F)\cong F_0$ where $F_0$ is the functor appearing in the defect sequence.  In fact, this isomorphism is natural in $F$ and hence we have the following.
\begin{prop*}$\extd_2\circ \extd_2\cong (\blank)_0$.\end{prop*}

We then show that the restriction of $\extd_2$ yields a duality of the defect zero functors.  At that point we derive the Generalised Auslander-Reiten formulas. We next use the localization properties of the functors $\agj$ and $\res$ to recover the duality $\newdual\cong \tr_*$.  We end by showing how to recover the Auslander-Reiten formulas from the Generalised Auslander-Reiten formulas.

\section{\label{section2}Functor Categories and The Defect}  

This section is a quick review.  For more details, see \cite{coh,stab} and \cite{mclane}.
Throughout this entire paper:
\begin{enumerate}
\item The word functor will always mean additive functor.    
\item The category of right modules will be denoted by $\Mod(R)$ and the category of finitely presented right modules will be denoted by $\mod(R)$.  
\item Every left $R$-module can be viewed as a right module over the opposite ring $R^{op}$.
\item The category of abelian groups will be denoted by $\ab$.
\item The category $\A$ will always denote an abelian category.
\end{enumerate}

\subsection{Finitely Presented Functors}

A functor $F\:\A\To \ab$ is called \dbf{representable} if it is isomorphic to $\hom_\A(X,\blank)$ for some $X\in \A$.   We will abbreviate the representable functors by $(X,\blank)$.   The most important property of representable functors is the following well known lemma of Yoneda:

\begin{lem}[Yoneda]For any covariant functor $F\:\A\to \ab$ and any $X\in \A$, there is an isomorphism:  $$\nat\big((X,\blank),F\big)\cong F(X)$$  given by $\alpha\mapsto \alpha_X(1_X)$.  The isomorphism is natural in both $F$ and $X$.  \end{lem}

An immediate consequence of the Yoneda lemma is that for any $X,Y\in \A$, $\nat\big((Y,\blank),(X,\blank)\big)\cong (X,Y)$.  Hence all natural transformations between representable functors come from morphisms between objects in $\A$.

\begin{defn}[Auslander]A functor $F\:\A\to \ab$ is called \textbf{finitely presented} if there exists a sequence of natural transformations \begin{center}\begin{tikzpicture}
\matrix(m)[ampersand replacement=\&, matrix of math nodes, row sep=3em, column sep=2em, text height=1.5ex, text depth=0.25ex]
{(Y,\blank)\&(X,\blank)\&F\&0\\}; 
\path[->,thick, font=\scriptsize]
(m-1-1) edge node[auto]{}(m-1-2)
(m-1-2) edge node[auto]{}(m-1-3)
(m-1-3) edge node[auto]{}(m-1-4);
\end{tikzpicture}\end{center}   such that for any $A\in \A$, the sequence of abelian groups \begin{center}\begin{tikzpicture}
\matrix(m)[ampersand replacement=\&, matrix of math nodes, row sep=3em, column sep=2em, text height=1.5ex, text depth=0.25ex]
{(Y,A)\&(X,A)\&F(A)\&0\\}; 
\path[->,thick, font=\scriptsize]
(m-1-1) edge node[auto]{}(m-1-2)
(m-1-2) edge node[auto]{}(m-1-3)
(m-1-3) edge node[auto]{}(m-1-4);
\end{tikzpicture}\end{center}  is exact.\end{defn}

One easily shows that if $F$ is finitely presented, then the collection of natural transformations $\nat(F,G)$ for any functor $G\:\A\To \ab$ is actually an abelian group.  As such, one may form a category whose objects are the covariant finitely presented functors $F\:\A\To \ab$ and whose morphisms are the natural transformations between two such functors.  This category is denoted by $\fp(\A,\ab)$ and was studied extensively by Auslander in multiple works.  

\begin{thm}[Auslander, \cite{coh}, Theorem 2.3]The category $\fp(\A,\ab)$ consisting of all finitely presented functors together with the natural transformations between them satisfies the following properties:
\begin{enumerate}
\item $\fp(\A,\ab)$ is abelian.  A sequence of finitely presented functors $$F\To G\To H$$ is exact if and only if for every $A\in \A$ the sequence of abelian groups $$F(A)\To G(A)\To H(A)$$ is exact.
\item The projective objects of $\fp(\A,\ab)$ are exactly the representables.   

\item Every finitely presented functor $F\in \fp(\A,\ab)$ has a projective resolution of the form $$0\To (Z,\blank)\To (Y,\blank)\To (X,\blank)\To F\To 0$$

\end{enumerate}
\end{thm}

The category $\fp(\A,\ab)$ also has enough injectives whenever $\A$ has enough projectives.  This was shown by Gentle and he completely classified them.

\begin{prop}[Gentle, \cite{rongentle}]\label{gentleinj}The category $\fp(\A,\ab)$ has enough injectives whenever $\A$ has enough projectives. A functor $F\:\A\To \ab$ is an injective object in $\fp(\A,\ab)$ if and only if there are projectives $P,Q\in \A$ such that there exists an exact sequence $$(Q,\blank)\to (P,\blank)\to F\to 0$$
\end{prop}   

\subsection{The Yoneda Embedding and the Defect}

The contravariant functor $\Y\:\A\To \fp(\A,\ab)$ given by $\Y(X):=(X,\blank)$ is a contravariant left exact embedding commonly referred to as the \dbf{Yoneda}\dbf{ embedding}.   At this point we will explicitly recall the construction of the \dbf{defect} functor $w\:\fp(\A,\ab)\to \A$, originally introduced by Auslander in \cite{coh}.  Let $F\in \fp(\A,\ab)$ and let  
\begin{center}\begin{tikzpicture}
\matrix(m)[ampersand replacement=\&, matrix of math nodes, row sep=3em, column sep=2.5em, text height=1.5ex, text depth=0.25ex]
{(Y,\blank)\&(X,\blank)\&F\&0\\}; 
\path[->,thick, font=\scriptsize]
(m-1-1) edge node[auto]{$(f,\blank)$}(m-1-2)
(m-1-2) edge node[auto]{$\alpha$}(m-1-3)
(m-1-3) edge node[auto]{}(m-1-4);
\end{tikzpicture}\end{center} be any presentation. The value of $w$ at $F$ is defined by the exact sequence 
\begin{center}\begin{tikzpicture}
\matrix(m)[ampersand replacement=\&, matrix of math nodes, row sep=3em, column sep=2em, text height=1.5ex, text depth=0.25ex]
{0\&w(F)\&X\&Y\\}; 
\path[->,thick, font=\scriptsize]
(m-1-1) edge node[auto]{}(m-1-2)
(m-1-2) edge node[auto]{$k$}(m-1-3)
(m-1-3) edge node[auto]{$f$}(m-1-4);
\end{tikzpicture}\end{center}This assignment extends to a contravariant exact functor \begin{center}\begin{tikzpicture}
\matrix(m)[ampersand replacement=\&, matrix of math nodes, row sep=3em, column sep=2em, text height=1.5ex, text depth=0.25ex]
{w\:\fp(\A,\ab)\&\A\\}; 
\path[->,thick, font=\scriptsize]
(m-1-1) edge node[auto]{}(m-1-2);
\end{tikzpicture}\end{center} To any projective resolution  \begin{center}\begin{tikzpicture}
\matrix(m)[ampersand replacement=\&, matrix of math nodes, row sep=1em, column sep=2.5em, text height=1.5ex, text depth=0.25ex]
{0\&(Z,\blank)\&(Y,\blank)\&(X,\blank)\&F\&0\\}; 
\path[->,thick, font=\scriptsize]
(m-1-1) edge node[auto]{}(m-1-2)
(m-1-2) edge node[auto]{$(g,\blank)$}(m-1-3)
(m-1-3) edge node[auto]{$(f,\blank)$}(m-1-4)
(m-1-4) edge node[auto]{$\alpha$}(m-1-5)
(m-1-5) edge node[auto]{}(m-1-6);
\end{tikzpicture}\end{center} apply the exact functor $w$.  This yields the exact sequence \begin{center}\begin{tikzpicture}
\matrix(m)[ampersand replacement=\&, matrix of math nodes, row sep=3em, column sep=2.5em, text height=1.5ex, text depth=0.25ex]
{0\&w(F)\&X\&Y\&Z\&0\\}; 
\path[->,thick, font=\scriptsize]
(m-1-1) edge node[auto]{}(m-1-2)
(m-1-2) edge node[auto]{$k$}(m-1-3)
(m-1-3) edge node[auto]{$f$}(m-1-4)
(m-1-4) edge node[auto]{$g$}(m-1-5)
(m-1-5) edge node[auto]{}(m-1-6);
\end{tikzpicture}\end{center}

\begin{diag}[Defect Diagram]This exact sequence embeds into the following commutative diagram, which we will from this point on refer to as the defect diagram, with exact rows and columns: 

\begin{center}\begin{tikzpicture}[framed]
\matrix(m)[ampersand replacement=\&, matrix of math nodes, row sep=1.5em, column sep=1.5em, text height=1.5ex, text depth=0.25ex]
{\&\&0\\
\&w(F)\&w(F)\\
0\&w(F)\&X\&Y\&Z\&0\\
\&0\&V\&Y\&Z\&0\\
\&\&0\\}; 
\path[->,thick, font=\scriptsize]
(m-1-3)edge (m-2-3)
(m-2-2) edge node[left]{$1$}(m-3-2)
(m-2-3) edge node[left]{$k$}(m-3-3)
(m-2-2) edge node[above]{$1$}(m-2-3)
(m-3-1) edge node[auto]{}(m-3-2)
(m-3-2) edge node[auto]{$k$}(m-3-3)
(m-3-3) edge node[auto]{$f$}(m-3-4)
(m-3-4) edge node[auto]{$g$}(m-3-5)
(m-3-5) edge node[auto]{}(m-3-6)
(m-4-2) edge node[auto]{}(m-4-3)
(m-4-3) edge node[auto]{$m$}(m-4-4)
(m-4-4) edge node[auto]{$g$}(m-4-5)
(m-4-5) edge node[auto]{}(m-4-6)
(m-3-3) edge node[left]{$e$}(m-4-3)
(m-3-4) edge node[left]{$1$}(m-4-4)
(m-3-5) edge node[left]{$1$}(m-4-5)
(m-4-3) edge (m-5-3);
\end{tikzpicture}\end{center}\end{diag}

Applying the Yoneda embedding to this diagram and extending to include cokernels where necessary yields the following commutative diagram with exact rows and columns: 

\begin{center}\begin{tikzpicture}
\matrix(m)[ampersand replacement=\&, matrix of math nodes, row sep=2em, column sep=2em, text height=1.5ex, text depth=0.25ex]
{\&\&\&0\&0\\
0\&(Z,\blank)\&(Y,\blank)\&(V,\blank)\&F_0\&0\\
0\&(Z,\blank)\&(Y,\blank)\&(X,\blank)\&F\&0\\
\&\&\&(w(F),\blank)\&(w(F),\blank)\\
\&\&\&F_1\&F_1\\
\&\&\&0\&0\\};
\path[->,thick,font=\scriptsize]
(m-1-4) edge (m-2-4)
(m-1-5) edge (m-2-5)
(m-2-1) edge (m-2-2)
(m-2-2) edge node[above]{$(g,\blank)$} (m-2-3)
edge node[left]{$1$}(m-3-2)
(m-2-3) edge node[above]{$(m,\blank)$} (m-2-4)
edge node[left]{$1$}(m-3-3)
(m-2-4) edge (m-2-5)
edge node[right]{$(e,\blank)$}(m-3-4)
(m-2-5) edge (m-2-6)
edge (m-3-5)
(m-3-1) edge (m-3-2)
(m-3-2) edge node[above]{$(g,\blank)$}(m-3-3)
(m-3-3) edge node[above]{$(f,\blank)$}(m-3-4)
(m-3-4) edge node[above]{$\alpha$}(m-3-5)
(m-3-5) edge (m-3-6)
(m-3-4) edge node[left]{$(k,\blank)$}(m-4-4)
(m-4-4) edge (m-5-4)
(m-5-4) edge (m-6-4)
(m-3-5) edge node[right]{$\varphi$}(m-4-5)
(m-4-5) edge (m-5-5)
(m-5-5) edge (m-6-5)
(m-4-4) edge node[above]{$1$}(m-4-5)
(m-5-4) edge node[above]{$1$}(m-5-5);
\end{tikzpicture}\end{center}

\begin{diag}[Defect Sequence]This yields the following exact sequence called the \dbf{defect sequence}: 
\begin{center}\begin{tikzpicture}
\matrix(m)[ampersand replacement=\&, matrix of math nodes, row sep=3em, column sep=2.5em, text height=1.5ex, text depth=0.25ex]
{0\&F_0\&F\&\big(w(F),\blank\big)\&F_1\&0\\}; 
\path[->,thick, font=\scriptsize]
(m-1-1) edge node[auto]{}(m-1-2)
(m-1-2) edge node[auto]{}(m-1-3)
(m-1-3) edge node[auto]{$\varphi$}(m-1-4)
(m-1-4) edge node[auto]{}(m-1-5)
(m-1-5) edge node[auto]{}(m-1-6);
\end{tikzpicture}\end{center}\end{diag}

\begin{prop}[Auslander, \cite{coh}, pg 204]\label{defadj}The Yoneda functor and the defect together form an adjoint pair \begin{center}\begin{tikzpicture}
\matrix (m) [ampersand replacement= \&,matrix of math nodes, row sep=3em, 
column sep=3em,text height=1.5ex,text depth=0.25ex] 
{\fp(\A,\ab)\&\A^{op}\\}; 
\path[->,thick, font=\scriptsize]
(m-1-1) edge[bend left=50] node[auto]{$w$}(m-1-2)
(m-1-2) edge[bend left=50] node[below]{$\Y$}(m-1-1);
\end{tikzpicture}\end{center}  The morphism $\varphi$ in the defect sequence is the unit of adjunction and $w(\varphi)$ is an isomorphism.  The functor $w$ is  exact and $w\Y\cong 1_\A$. \end{prop}

Because the defect $w$ is an exact functor, its kernel, that is the full subcategory of $\fp(\A,\ab)$ consisting of all functors $F$ for which $w(F)=0$, is a Serre subcategory of $\fp(\A,\ab)$.  We will denote this category by both $\ker(w)$ and by $\fp_0(\A,\ab)$ depending on the situation.  Auslander studied $\fp_0(\A,\ab)$ in \cite{coh} and established a list of properties which we will be useful for certain portions of our discussion.

\begin{lem}[Auslander, \cite{coh}, Lemma 4.1, Lemma 4.2]If $\A$ has enough projectives, then the following properties hold.
\begin{enumerate}
\item All extension functors $\ext^1(A,\blank)$ are in $\fp(\A,\ab)$ and in fact in $\fp_0(\A,\ab)$.
\item The functors $\ext^1(A,\blank)$ are injective objectives in $\fp_0(\A,\ab)$.
\item For every functor $F\in \fp_0(\A,\ab)$, there is a monomorphism $$F\to \ext^1(A,\blank)$$ for some $A\in \A$.  That is the Serre subcategory $\fp_0(\A,\ab)$ has enough injectives of the form $\ext^1(A,\blank)$.
\item A functor $F$ in $\fp_0(\A,\ab)$ is injective if and only if $F$ is half exact.
\end{enumerate}
\end{lem}

Moreover, Freyd showed in \cite{freyd} that if $\A$ has countable sums then any direct summand $F$ of an extension functor $\ext^1(A,\blank)$ is isomorphic to $\ext^1(C,\blank)$ for some $C\in \A$.  Combining this with Auslander's result, if $\A$ is has enough projectives and has countable sums, then any half exact finitely presented functor $F\:\A\To \ab$ is isomorphic to an extension functor $\ext^1(C,\blank)$ for some $C\in \A$.  Thus in the case $\A=\Mod(R)$, the category $\fp_0(\A,\ab)$ has enough injectives and they are precisely the extension functors.

\subsection{The Tensor Embedding}

The Yoneda lemma allows one to study the functor category $\fp(\A,\ab)$ for any abelian category $\A$ due to the fact that the natural transformations between any two such functors form an actual set.  For any ring $R$, the module category $\Mod(R)$ can be embedded into the functor category $\fp(\Mod(R),\ab)$ using the Yoneda embedding.  This allows one to study the category $\Mod(R)$ by looking at the category of finitely presented functors.    There are other approaches to studying the category $\Mod(R)$ using a suitable functor category.

Denote by $(\mod(R^{op}),\ab)$ the category whose objects are additive functors $F\:\mod(R^{op})\To \ab$ and whose morphisms are the natural transformations between.   This category is also abelian and again a sequence of functors $F\to G\to H$ in $(\mod(R^{op}),\ab)$ is exact if and only if for any finitely presented module $A\in \mod(R^{op})$, the sequence of abelian groups $F(A)\to G(A)\to H(A)$ is exact.

The category $(\mod(R^{op}),\ab)$ is in fact Grothendieck. Hence it has enough projectives and injectives.   The projectives of $(\mod(R^{op}),\ab)$ are the summands of functors of the form $\bigoplus_{i\in I}(X_i,\blank)$ such that each $X_i\in \mod(R^{op})$.  The injectives were classified by Gruson and Jensen in \cite{gj} and we now recall their classification.  

A  \dbf{pure exact sequence} in $\Mod(R)$ is a short exact sequence of modules $$0\to A\to B\to C\to 0$$ such that $$0\to A\otimes\blank\to B\otimes\blank\to C\otimes\blank\to 0$$ is exact in $(\mod(R^{op}),\ab)$.  A \dbf{pure injective} module $M$ is any module which is injective with respect to pure exact sequences.    That is given any pure exact sequence  $$0\to A\to B\to C\to 0$$ the sequence $$0\to(C,M)\to (B,M)\to(A,M)\to 0$$ is exact.

\begin{thm}[Gruson, Jensen, \cite{gj}]The functor $t\:\Mod(R)\To (\mod(R^{op}),\ab)$ given by $$t(M)=M\otimes\blank$$ satisfies the following two properties.
\begin{enumerate}
\item  $t$ is a fully faithful covariant functor.  Due to this property, we call the functor $t$ the \dbf{tensor embedding}.   \item The injectives of $(\mod(R^{op}),\ab)$ are precisely the functors of the form $M\otimes\blank$ for which $M$ is pure injective.  That is, a functor $F\:\mod(R^{op})\to \ab$ is injective in $(\mod(R^{op}),\ab)$ if and only if it is isomorphic to $t(M)$ for some pure injective $M\in\Mod(R)$.
\end{enumerate}  \end{thm}

\subsection{The Auslander-Gruson-Jensen Functor}

In \cite{dr1} we study the relationship between the functor categories $\fp(\Mod(R),\ab)$ and  $(\mod(R^{op}),\ab)$.  
The category $(\mod(R^{op}),\ab)$ is abelian and the category $\fp(\Mod(R),\ab)$ has enough projectives.  Given any contravariant functor $$S\:\fp(\Mod(R),\ab)\To (\mod(R^{op}),\ab)$$ we can use the projectives to compute the zeroth right derived functor $R_0S$ which will be a contravariant left exact functor.   By composing the tensor embedding $t$ with the defect $w$ we get a functor $$t\circ w\:\fp(\Mod(R),\ab)\To (\mod(R^{op}),\ab)$$ which is right exact contravariant and sends any finitely presented functor $F$ to the tensor functor $w(F)\otimes\blank$.  

The \dbf{Auslander-Gruson-Jensen functor}\footnote{The functor $\agj$ is named as such due to the fact that Auslander in \cite{isosing} and independently Gruson and Jensen in \cite{gj} define for any additive functor $F\:\Mod(R)\to \ab$, a corresponding functor $DF\:\mod(R^{op})\to \ab$ which agrees with $\agj F$ in the case that $F$ is finitely presented.} is the zeroth right derived functor of $t\circ w$. That is,  
$$\agj=R_0(t\circ w)$$  As shown in \cite{dr1}, the Auslander-Gruson-Jensen functor $$\agj\:\fp(\Mod(R),\ab)\to (\mod(R^{op}),\ab)$$ is completely determined by  the following
\begin{enumerate}
\item $\agj$ is exact contravariant.
\item For any $(X,\blank)$, $\agj(X,\blank)= X\otimes\blank$. 
\item For any $(f,\blank)\:(Y,\blank)\To (X,\blank)$, $\agj(f,\blank)=f\otimes\blank$.
%\item For any $L\in \mod(R^{op})$, $\agj(\blank\otimes L)\cong (L,\blank)$.
\end{enumerate}

\subsection{Stable Categories}Suppose that $\A$ is an additive category with enough projectives.  Given any two objects $X,Y\in \A$ the set $\P(X,Y)$ is defined to be the collection of all $f\in \hom_\A(X,Y)$ such that there exists a projective $P$ and a commutative triangle 

\begin{center}\begin{tikzpicture}
\matrix(m)[ampersand replacement=\&, matrix of math nodes, row sep=1.5em, column sep=1em, text height=0.5ex, text depth=0.25ex]
{X\&\&Y\\
\&P\\};
\path[->,thick, font=\scriptsize]
(m-1-1) edge node[above]{$f$} (m-1-3)
(m-1-1) edge node[left]{$a$} (m-2-2)
(m-2-2) edge node[right]{$b$}(m-1-3);
\end{tikzpicture}\end{center} $\P(X,Y)$ is a subgroup of $\hom(X,Y)$ and hence one can move to the quotient group $$\underline{\hom}(X,Y)=\frac{\hom(X,Y)}{\P(X,Y)}$$  The \dbf{projectively}  \dbf{stable category} $\underline{\A}$ is defined as follows.  The objects of  $\underline{\A}$ are the same as that of the category $\A$.  Given $X,Y\in \A$, the morphism set $\hom_{\underline{\A}}(X,Y)$ is defined as follows $$ \hom_{\underline{\A}}(X,Y)=\underline{\hom}(X,Y)$$  

Dually one can create for any additive category with enough injectives, the \dbf{injectively stable category} $\overline{\A}$. The category $\overline{\A}$ comes with a functor $Q\:\A\to \overline{\A}$ sending every object $A\in \A$ to the object $A\in \overline{\A}$.  The functor $Q$ together with $\overline{\A}$ satisfy the following universal property.  If $\B$ is any additive category and $\E\:\A\to \B$ is any additive functor vanishing on all injectives in $\A$, then there is a unique $\Psi$ such that the following diagram commutes

\begin{center}\begin{tikzpicture}
\matrix(m)[ampersand replacement=\&, matrix of math nodes, row sep=3em, column sep=3em, text height=1.5ex, text depth=0.25ex]
{\A\&\overline{\A}\\
\B\&\\};
\path[->,thick, font=\scriptsize]
(m-1-1) edge node[left ]{$ \E$} (m-2-1)
(m-1-2) edge node[auto]{$\exists !\Psi$} (m-2-1)
(m-1-1) edge node[auto]{$Q$}(m-1-2);
\end{tikzpicture}\end{center}  The category $\underline{\A}$ also comes with a functor $Q\:\A\to \underline{\A}$ satisfying a similar universal property.

\section{\label{section3}Reflective Adjunctions and the Derived Recollement Theorem}
\textbf{All categories considered in this section will be abelian.}  Some of the statements may hold in generality; however, our focus is on the abelian setting.  The reader is referred to Theorem 1 in Chapter IV, Section 3 from \cite{mclane}, for a proof of the following lemma, which we will use without comment.

\begin{lem}\label{ffad}If $(\la,\ra,\eta,\epsilon)$ is an adjunction \begin{center}\begin{tikzpicture}
\matrix (m) [ampersand replacement= \&,matrix of math nodes, row sep=3em, 
column sep=4em,text height=1.5ex,text depth=0.25ex] 
{\C\&\D\\}; 
\path[->,thick, font=\scriptsize]
(m-1-1) edge[bend left=50] node[auto]{$\la$}(m-1-2)
(m-1-2) edge[bend left=50] node[below]{$\ra$}(m-1-1);
\end{tikzpicture}\end{center}with unit $\eta$ and counit $\eps$, then 
\begin{enumerate}
%\item $\ra$ is full if and only if every component $\eps_A$ is an epimorphism.
%\item $\ra$ is faithful if and only if every component $\eps_A$ is a split monomorphism.
\item $\ra$ is fully faithful if and only if the counit of adjunction $\eps$ is an isomorphism.
\item $\la$ is fully faithful if and only if the unit of adjunction $\eta$ is an isomorphism.
\end{enumerate}
\end{lem}

\begin{defn}
A \dbf{recollement} consists of a diagram of functors between abelian categories $\A,\B,\C$,

\begin{center}\begin{tikzpicture}
\matrix (m) [ampersand replacement= \&,matrix of math nodes, row sep=5em, 
column sep=4em,text height=1.5ex,text depth=0.25ex] 
{\A\&\B\&\C\\}; 
\path[->,thick, font=\scriptsize]
(m-1-1) edge node[above]{$i$}(m-1-2)
(m-1-2) edge[bend right=50] node[above]{$i_L$}(m-1-1)
(m-1-2) edge[bend left=50] node[below]{$i_R$}(m-1-1)
(m-1-3)edge[bend left=50] node[below]{$f_R$}(m-1-2)
(m-1-3) edge[bend right=50] node[above]{$f_L$}(m-1-2)
(m-1-2) edge node[above]{$f$}(m-1-3);
\end{tikzpicture}\end{center} where 
\begin{enumerate}
\item The functor $f$ has left adjoint $f_L$ and right adjoint $f_R$.
\item The unit $1_\C\to ff_L$ and counit $ff_R\to 1_\C$ are isomorphisms.
\item The functor $i$ has left adjoint $i_L$ and right adjoint $i^R$.
\item The unit $1_\A\to i_Ri$ and counit $i_Li\to 1_\A$ are isomorphisms.
\item The functor $i$ is an embedding onto the full subcategory of $\A$ consisting of all objects $A$ for which $f(A)=0$.
\end{enumerate}
Equivalently, a recollement is such a diagram with:
\begin{enumerate}
\item The functor $f$ has left adjoint $f_L$ and right adjoint $f_R$.
\item $f_L$ and $f_R$ are fully faithful.
\item The functor $i$ has left adjoint $i_L$ and right adjoint $i_R$.
\item $i$ is fully faithful.
\item The essential image of $i$, $\{B\in \mathcal{B}:B\cong i(A)\text{ for some }A\in\mathcal{A}\}$, is equal to the kernel of $f$, $\{B\in\mathcal{B}:f(B)=0\}$.
\end{enumerate}
\end{defn}

\begin{defn}An adjunction $\la\dashv\ra$
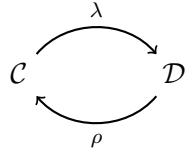
\begin{center}\begin{tikzpicture}
\matrix (m) [ampersand replacement= \&,matrix of math nodes, row sep=3em, 
column sep=4em,text height=1.5ex,text depth=0.25ex] 
{\C\&\D\\}; 
\path[->,thick, font=\scriptsize]
(m-1-1) edge[bend left=50] node[auto]{$\la$}(m-1-2)
(m-1-2) edge[bend left=50] node[below]{$\ra$}(m-1-1);
\end{tikzpicture}\end{center}  is called \dbf{reflective} if $\ra$ is fully faithful and \dbf{coreflective} if $\la$ is fully faithful. Equivalently, it is reflective iff the counit $\la\ra\to 1$ is an isomorphism, and it is coreflective iff the unit $1\to\ra\la$ is an isomorphism.
\end{defn}

\begin{lem}\label{samprop}Given a functor $G:\C\to \D$ with adjoints $F\dashv G\dashv H$, $F$ is fully faithful if and only if $H$ is fully faithful.  In the language above, the adjunction $F\dashv G$ is coreflective if and only if $G\dashv H$ is reflective.\end{lem}

\begin{pf}Let $\eta$ and $\epsilon$ be the unit and counit of the adjunction $F\dashv G$ and let $\theta$ and $\zeta$ be the unit and counit of the adjunction $G\dashv H$.

If $F$ is fully faithful, the unit $\eta:1_\D\to GF$ is an isomorphism. Therefore $GF$ is an equivalence. By composition of adjunctions, it follows that $GF\dashv GH$, with unit $(G\theta F)\eta$ and counit $\zeta(H\epsilon G)$. Since $GF$ is an equivalence, so is $GH$, and hence both are fully faithful. Therefore, both $(G\theta F)\eta$ and $\zeta(H\epsilon G)$ are isomorphisms. Note that, since $(\eta G)(\epsilon G)=1_G$, $\epsilon G$ and hence $H\epsilon G$ is an isomorphism. It follows that $\zeta$ is an isomorphism, and hence $H$ is fully faithful.  The other direction is dual.$\qed$ \end{pf}

The following generalizes Auslander's approach to constructing the defect sequence from \cite{coh}.

\begin{lem}Let  $(\la,\ra,\eta,\eps)$ denote an adjunction
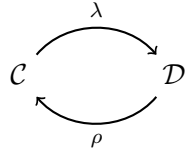
\begin{center}\begin{tikzpicture}
\matrix (m) [ampersand replacement= \&,matrix of math nodes, row sep=3em, 
column sep=4em,text height=1.5ex,text depth=0.25ex] 
{\C\&\D\\}; 
\path[->,thick, font=\scriptsize]
(m-1-1) edge[bend left=50] node[auto]{$\la$}(m-1-2)
(m-1-2) edge[bend left=50] node[below]{$\ra$}(m-1-1);
\end{tikzpicture}\end{center} 
If this adjunction is reflective, then  we have the following.
\begin{enumerate}
\item $\la(\eta_C)$ is an isomorphism for any $C$.
\item For each $C\in \C$, there exists an exact sequence natural in $C$
\begin{center}\begin{tikzpicture}
\matrix(m)[ampersand replacement=\&, matrix of math nodes, row sep=2em, column sep=2em, text height=1.5ex, text depth=0.25ex]
{0\&C_0\&C\&\rho\circ\lambda(C)\&C_1\&0\\}; 
\path[->,thick, font=\scriptsize]
(m-1-1) edge node[auto]{}(m-1-2)
(m-1-2) edge node[auto]{}(m-1-3)
(m-1-3) edge node[auto]{$\eta_C$}(m-1-4)
(m-1-5) edge node[auto]{}(m-1-6)
(m-1-4) edge node[auto]{}(m-1-5);
\end{tikzpicture}\end{center} such that $\la(C_1)\cong 0$.  
\item If $\la$ is exact, then $\la(C_0)\cong 0$ and the functor $(\blank)_0\:\C\to \ker(\la)$ is the right adjoint to the inclusion functor $i\:\ker(\la)\to \C$.  That is there is an adjoint pair

\begin{center}\begin{tikzpicture}
\matrix (m) [ampersand replacement= \&,matrix of math nodes, row sep=3em, 
column sep=4em,text height=1.5ex,text depth=0.25ex] 
{\ker(\la)\&\C\\}; 
\path[->,thick, font=\scriptsize]
(m-1-1) edge[bend left=50] node[auto]{$i$}(m-1-2)
(m-1-2) edge[bend left=50] node[below]{$(\blank)_0$}(m-1-1);
\end{tikzpicture}\end{center}
\end{enumerate}
\end{lem}

\begin{pf}By the triangle identity for adjunctions the following composition for any $C\in \C$ is the identity.
 \begin{center}\begin{tikzpicture}
\matrix (m) [ampersand replacement= \&,matrix of math nodes, row sep=3em, 
column sep=4em,text height=1.5ex,text depth=0.25ex] 
{\la(C)\&\la\circ\ra\circ\la(C)\&\la(C)\\}; 
\path[->,thick, font=\scriptsize]
(m-1-1) edge node[above]{$\la(\eta_C)$}(m-1-2)
(m-1-2) edge node[above]{$\eps_{\la(C)}$}(m-1-3);
\end{tikzpicture}\end{center} Since $\ra$ is fully faithful, $\eps_{\la(C)}$ is an isomorphism.  Therefore, $\la(\eta_C)$ is also an isomorphism. 

The exact sequence
\begin{center}\begin{tikzpicture}
\matrix(m)[ampersand replacement=\&, matrix of math nodes, row sep=2em, column sep=2em, text height=1.5ex, text depth=0.25ex]
{0\&C_0\&C\&\ra\circ\la(C)\&C_1\&0\\}; 
\path[->,thick, font=\scriptsize]
(m-1-1) edge node[auto]{}(m-1-2)
(m-1-2) edge node[auto]{}(m-1-3)
(m-1-3) edge node[auto]{$\eta_C$}(m-1-4)
(m-1-5) edge node[auto]{}(m-1-6)
(m-1-4) edge node[auto]{}(m-1-5);
\end{tikzpicture}\end{center} is obtained by taking the kernel and cokernel of $\eta_C$.  The fact that $\la(\eta_C)$ is an isomorphism and that $\la$ is always right exact forces $\la(C_1)\cong 0$.  

If $\la$ is exact, then $\la(C_0)\cong 0$ again because $\la(\eta_C)$ is an isomorphism.  It is also clear that if $\la(B)\cong 0$, then $B\cong B_0$.  Therefore, for any object  $B\in \ker(\la)$, there is an isomorphism $$\hom(i(B),C)\cong\hom_{\ker(\la)}(B,C_0)$$ which is easily seen to be natural in $B$ and $C$.  $\qed$
\end{pf}

We state the dual statement without proof.

\begin{lem}Let  $(\lla,\la,\eta,\eps)$ denote an adjunction\begin{center}\begin{tikzpicture}
\matrix (m) [ampersand replacement= \&,matrix of math nodes, row sep=3em, 
column sep=4em,text height=1.5ex,text depth=0.25ex] 
{\B\&\C\\}; 
\path[->,thick, font=\scriptsize]
(m-1-1) edge[bend left=50] node[auto]{$\lla$}(m-1-2)
(m-1-2) edge[bend left=50] node[below]{$\la$}(m-1-1);
\end{tikzpicture}\end{center} 
If this adjunction is coreflective, then  we have the following.
\begin{enumerate}
\item $\la(\eps_C)$ is an isomorphism for any $C$.
\item For each $C\in \C$, there exists an exact sequence natural in $C$
\begin{center}\begin{tikzpicture}
\matrix(m)[ampersand replacement=\&, matrix of math nodes, row sep=2em, column sep=2em, text height=1.5ex, text depth=0.25ex]
{0\&C^1\&\lla\circ \la(C)\&C\&C^0\&0\\}; 
\path[->,thick, font=\scriptsize]
(m-1-1) edge node[auto]{}(m-1-2)
(m-1-2) edge node[auto]{}(m-1-3)
(m-1-3) edge node[auto]{$\eps_C$}(m-1-4)
(m-1-5) edge node[auto]{}(m-1-6)
(m-1-4) edge node[auto]{}(m-1-5);
\end{tikzpicture}\end{center} such that $\la(C^1)\cong 0$.  
\item If $\la$ is exact, then $\la(C^0)\cong 0$ and the functor $(\blank)^0\:\C\to \ker(\la)$ is the left adjoint to the inclusion functor $i\:\ker(\la)\to \C$.  That is there is an adjoint pair

\begin{center}\begin{tikzpicture}
\matrix (m) [ampersand replacement= \&,matrix of math nodes, row sep=3em, 
column sep=4em,text height=1.5ex,text depth=0.25ex] 
{\ker(\la)\&\C\\}; 
\path[->,thick, font=\scriptsize]
(m-1-1) edge[bend right=50] node[below]{$i$}(m-1-2)
(m-1-2) edge[bend right=50] node[above]{$(\blank)^0$}(m-1-1);
\end{tikzpicture}\end{center}
\end{enumerate}
\end{lem}

\begin{cor}\label{abcritrec}Suppose that we have adjunctions $(\la,\ra, \eta,\eps)$ and $(\lla,\la,\overline{\eta},\overline{\eps})$

\begin{center}\begin{tikzpicture}
\matrix (m) [ampersand replacement= \&,matrix of math nodes, row sep=3em, 
column sep=4em,text height=1.5ex,text depth=0.25ex] 
{\C\&\D\\}; 
\path[->,thick, font=\scriptsize]
(m-1-2) edge[bend right=50]  node[above]{$\lla$}(m-1-1)
(m-1-1) edge node[above]{$\la$}(m-1-2)
(m-1-2) edge[bend left=50] node[below]{$\ra$}(m-1-1);
\end{tikzpicture}\end{center} The following are equivalent
\begin{enumerate}
\item $(\lla,\la,\overline{\eta},\overline{\eps})$ is coreflective.
\item $(\la,\ra, \eta,\eps)$ is reflective. 
\item We have a recollement 

\begin{center}\begin{tikzpicture}
\matrix (m) [ampersand replacement= \&,matrix of math nodes, row sep=5em, 
column sep=5em,text height=1.5ex,text depth=0.25ex] 
{\ker(\la)\&\C\&\D\\}; 
\path[->,thick, font=\scriptsize]
(m-1-1) edge node[above]{$i$}(m-1-2)
(m-1-2) edge[bend right=50] node[above]{$(\blank)^0$}(m-1-1)
(m-1-2) edge[bend left=50] node[below]{$(\blank)_0$}(m-1-1)
(m-1-3)edge[bend left=50] node[below]{$\ra$}(m-1-2)
(m-1-3) edge[bend right=50] node[above]{$\lla$}(m-1-2)
(m-1-2) edge node[above]{$\la$}(m-1-3);
\end{tikzpicture}\end{center} 
\end{enumerate}
\end{cor}

\begin{pf}(3) clearly implies (1) and (2).  (1) and (2) are equivalent by Lemma \ref{samprop}.  Since $\la$ is a left and right adjoint, it is exact.  Now apply the preceding lemmas.  $\qed$\end{pf}

The following is well known but is included for the purposes of self-containment.  

\begin{prop}\label{serloc}If the adjunction $(\la,\ra,\eta,\eps)$ \begin{center}\begin{tikzpicture}
\matrix (m) [ampersand replacement= \&,matrix of math nodes, row sep=3em, 
column sep=4em,text height=1.5ex,text depth=0.25ex] 
{\C\&\D\\}; 
\path[->,thick, font=\scriptsize]
(m-1-1) edge[bend left=50] node[auto]{$\la$}(m-1-2)
(m-1-2) edge[bend left=50] node[below]{$\ra$}(m-1-1);
\end{tikzpicture}\end{center} is reflective and $\la$ is exact, then $\la$ is a Serre localization functor in the following sense. 
For every exact functor $\E\:\C\To \A$ into abelian category $\A$ such that $\ker(\E)\sub \ker(\la)$, there exists a unique functor $\psi=\E\ra$ making the following diagram commute. Furthermore, $\psi$ is exact.

\begin{center}\begin{tikzpicture}
\matrix(m)[ampersand replacement=\&, matrix of math nodes, row sep=3em, column sep=3em, text height=1.5ex, text depth=0.25ex]
{\C\&\D\\
\B\&\\};
\path[->,thick, font=\scriptsize]
(m-1-1) edge node[left ]{$ \E$} (m-2-1)
(m-1-2) edge node[auto]{$\exists !\psi$} (m-2-1)
(m-1-1) edge node[auto]{$\la$}(m-1-2);
\end{tikzpicture}\end{center}  This induces an equivalence of categories $$\frac{\C}{\ker{\la}}\cong \D$$

\end{prop}

\begin{pf}Define $\psi=\E\ra$.  We begin by showing that $\psi$ is exact.  Take any short exact sequence in $\D$

\begin{center}\begin{tikzpicture}
\matrix(m)[ampersand replacement=\&, matrix of math nodes, row sep=2em, column sep=2em, text height=1.5ex, text depth=0.25ex]
{0\&X\&Y\&Z\&0\\}; 
\path[->,thick, font=\scriptsize]
(m-1-1) edge node[auto]{}(m-1-2)
(m-1-2) edge node[auto]{}(m-1-3)
(m-1-3) edge node[auto]{}(m-1-4)
(m-1-4) edge node{}(m-1-5);
\end{tikzpicture}\end{center} Applying the right adjoint $\ra$ which is left exact functor $\ra$ yields the exact sequence in $\C$

\begin{center}\begin{tikzpicture}
\matrix(m)[ampersand replacement=\&, matrix of math nodes, row sep=2em, column sep=2em, text height=1.5ex, text depth=0.25ex]
{0\&\ra(X)\&\ra(Y)\&\ra(Z)\&C\&0\\}; 
\path[->,thick, font=\scriptsize]
(m-1-1) edge node[auto]{}(m-1-2)
(m-1-2) edge node[auto]{}(m-1-3)
(m-1-3) edge node[auto]{}(m-1-4)
(m-1-4) edge node{}(m-1-5)
(m-1-5) edge node{}(m-1-6);
\end{tikzpicture}\end{center}

Since this sequence is exact and $\la\ra\cong 1_\D$, applying the exact functor $\la$ returns us to the original short exact sequence \begin{center}\begin{tikzpicture}
\matrix(m)[ampersand replacement=\&, matrix of math nodes, row sep=2em, column sep=2em, text height=1.5ex, text depth=0.25ex]
{0\&X\&Y\&Z\&0\\}; 
\path[->,thick, font=\scriptsize]
(m-1-1) edge node[auto]{}(m-1-2)
(m-1-2) edge node[auto]{}(m-1-3)
(m-1-3) edge node[auto]{}(m-1-4)
(m-1-4) edge node{}(m-1-5);
\end{tikzpicture}\end{center} which shows that $\la(C)\cong 0$.  By assumption $\ker(\la)\sub \ker(\E)$.  Therefore $\E(C)\cong 0$.  Hence by applying $\E$ to the exact sequence
\begin{center}\begin{tikzpicture}
\matrix(m)[ampersand replacement=\&, matrix of math nodes, row sep=2em, column sep=2em, text height=1.5ex, text depth=0.25ex]
{0\&\ra(X)\&\ra(Y)\&\ra(Z)\&C\&0\\}; 
\path[->,thick, font=\scriptsize]
(m-1-1) edge node[auto]{}(m-1-2)
(m-1-2) edge node[auto]{}(m-1-3)
(m-1-3) edge node[auto]{}(m-1-4)
(m-1-4) edge node{}(m-1-5)
(m-1-5) edge node{}(m-1-6);
\end{tikzpicture}\end{center} one has the exact sequence

\begin{center}\begin{tikzpicture}
\matrix(m)[ampersand replacement=\&, matrix of math nodes, row sep=2em, column sep=2em, text height=1.5ex, text depth=0.25ex]
{0\&\E\ra(X)\&\E\ra(Y)\&\E\ra(Z)\&0\\}; 
\path[->,thick, font=\scriptsize]
(m-1-1) edge node[auto]{}(m-1-2)
(m-1-2) edge node[auto]{}(m-1-3)
(m-1-3) edge node[auto]{}(m-1-4)
(m-1-4) edge node{}(m-1-5);
\end{tikzpicture}\end{center}  This shows that $\psi=\E\ra$ is exact.

We want to show that $\psi$ is the unique exact functor that satisfies the relation $\psi \la\cong \E$.  First we show that it satisfies this relation. Take any object $C\in \C$ and apply the functor $\E$ to the exact sequence \begin{center}\begin{tikzpicture}
\matrix(m)[ampersand replacement=\&, matrix of math nodes, row sep=3em, column sep=2.5em, text height=1.5ex, text depth=0.25ex]
{0\&C_0\&C\&\ra\la(C)\&C_1\&0\\}; 
\path[->,thick, font=\scriptsize]
(m-1-1) edge node[auto]{}(m-1-2)
(m-1-2) edge node[auto]{}(m-1-3)
(m-1-3) edge node[auto]{$\eta_C$}(m-1-4)
(m-1-4) edge node[auto]{}(m-1-5)
(m-1-5) edge node[auto]{}(m-1-6);
\end{tikzpicture}\end{center}  To get the exact sequence 

\begin{center}\begin{tikzpicture}
\matrix(m)[ampersand replacement=\&, matrix of math nodes, row sep=3em, column sep=2.5em, text height=1.5ex, text depth=0.25ex]
{0\&\E(C_0)\&\E (C)\&\E(\ra\la(C))\&\E(C_1)\&0\\}; 
\path[->,thick, font=\scriptsize]
(m-1-1) edge node[auto]{}(m-1-2)
(m-1-2) edge node[auto]{}(m-1-3)
(m-1-3) edge node[auto]{$\E(\eta_C)$}(m-1-4)
(m-1-4) edge node[auto]{}(m-1-5)
(m-1-5) edge node[auto]{}(m-1-6);
\end{tikzpicture}\end{center}  Since $C_0,C_1\in \ker(\la)$, by assumption on $\E$ we also know that $\E(C_0)\cong \E(C_1)\cong 0$ which makes $\E(\eta_C)$ an isomorphism.  Since $\eta$ is natural in $F$, one now has that $\E \la \ra\cong \E$ establishing that $\psi \ra\cong \E$ as required.  

We conclude with the uniqueness portion of the argument.  If $\theta\:\D\To \A$ is exact and satisfies $\theta \ra\cong \E$, then since $\la\ra\cong 1_\D$  we have $$\theta\cong \theta 1_\D\cong \theta \la \ra\cong \E\ra=\psi$$ which shows that $\theta \cong \psi$ as claimed.  This completes the proof that $\la$ is a Serre localization functor.  $\qed$ \end{pf}

The next lemma will be used implicitly in our main theorem.
\begin{lem}Let $\alpha:F\to G$ be a natural transformation between right exact functors $F,G:\C\to\D$ where $\C$ and $\D$ are abelian categories and $\C$ has enough projectives. If $\alpha_P:F(P)\to G(P)$ is an isomorphism for any projective $P\in\C$ then $\alpha$ is an isomorphism.
\end{lem}

\begin{thm}[Derived Recollement Theorem]\label{recthm}Suppose that $\C,\D$ are both abelian categories with enough projectives and we have a reflective adjunction

\begin{center}\begin{tikzpicture}
\matrix (m) [ampersand replacement= \&,matrix of math nodes, row sep=3em, 
column sep=4em,text height=1.5ex,text depth=0.25ex] 
{\C\&\D\\}; 
\path[->,thick, font=\scriptsize]
(m-1-1) edge[bend left=50] node[auto]{$\la$}(m-1-2)
(m-1-2) edge[bend left=50] node[below]{$\ra$}(m-1-1);
\end{tikzpicture}\end{center} Assume that 
\begin{enumerate}
\item For any projective $P\in \C$ there is some object $D\in \D$ such that $P\cong\ra (D)$.
\item $\rho$ sends projectives to projectives.
\end{enumerate}

Then the adjunction $\la\dashv\ra$ can be extended to a recollement

\begin{center}\begin{tikzpicture}
\matrix (m) [ampersand replacement= \&,matrix of math nodes, row sep=5em, 
column sep=5em,text height=1.5ex,text depth=0.25ex] 
{\ker(\la)\&\C\&\D\\}; 
\path[->,thick, font=\scriptsize]
(m-1-1) edge node[above]{$i$}(m-1-2)
(m-1-2) edge[bend right=50] node[above]{$i_L$}(m-1-1)
(m-1-2) edge[bend left=50] node[below]{$i_R$}(m-1-1)
(m-1-3)edge[bend left=50] node[below]{$\ra$}(m-1-2)
(m-1-3) edge[bend right=50] node[above]{$L_0\ra$}(m-1-2)
(m-1-2) edge node[above]{$\la$}(m-1-3);
\end{tikzpicture}\end{center} 
In particular, $\la$ is exact.
\end{thm}

\begin{pf}By Corollary \ref{abcritrec} all that we must do is show that we have an  adjunction $(L_0\ra, \la)$.   Let $\mu=L_0\rho$. Let $\epsilon:\lambda\rho\to 1$ be the counit of the adjunction $(\lambda,\rho)$, which is an isomorphism. Since $\rho$ preserves projectives, we may assume $\mu (P)=\rho (P)$ for any projective $P\in \D$,

We first construct the counit $\zeta:1\to \lambda\mu$. Let $D\in \D$ be given, and choose a projective presentation $P\to Q\to D\to 0$. We define $\zeta_D:D\to\lambda\mu (D)$ to be the unique morphism such that the diagram
\begin{displaymath}
\xymatrix{P\ar[d]_{\epsilon^{-1}_P}\ar[r]&Q\ar[d]_{\epsilon^{-1}_Q}\ar[r]&D\ar[d]^{\zeta_D}\ar[r]&0\\
\lambda\rho (P)\ar[r]&\lambda\rho (Q)\ar[r]&\lambda\mu (D)\ar[r]&0
}
\end{displaymath}
commutes. Using the comparison theorem of homological algebra, one can show that $\zeta_D$ is natural in $D$. 

For any $D\in \D$ and $C\in \C$, we have the morphism
\begin{displaymath}
\hom(\mu (D),C)\to \hom(D,\la (C)):g\mapsto (\la (g))\zeta_D.
\end{displaymath}
We need to show that this is bijective. Since $\mu$ is right exact and $\D$ has enough projectives, we may assume that $D$ is projective. Therefore, $\mu (D)=\rho (D)$ and we need only show that, when $D$ is projective,
\begin{displaymath}
\hom(\ra (D),C)\to \hom(D,\la (C)):g\mapsto (\la (g))\zeta_D.
\end{displaymath}
But since, $\ra (D)$ is projective and $\la$ is right exact and $\C$ has enough projectives, we may assume that $C$ is projective, and hence $C=\rho (D')$ for an object $D'\in\D$. Since $\ra$ is fully faithful, our morphism becomes
\begin{align*}
\hom(\ra (D),\ra (D'))&\to \hom(D,\la (D'))\\
\ra (f)&\mapsto (\la\ra (f))\epsilon^{-1}_D=\epsilon^{-1}_{D'}f
\end{align*}
which, since $\ra$ is fully faithful, is clearly a bijection.

\end{pf}
We state without proof the dual.

\begin{thm}[Derived Recollement Theorem*]\label{retthmdual}Suppose that $\C,\D$ are both abelian categories with enough injectives and we have a coreflective adjunction

\begin{center}\begin{tikzpicture}
\matrix (m) [ampersand replacement= \&,matrix of math nodes, row sep=3em, 
column sep=4em,text height=1.5ex,text depth=0.25ex] 
{\C\&\D\\}; 
\path[->,thick, font=\scriptsize]
(m-1-1) edge[bend left=50] node[auto]{$\la$}(m-1-2)
(m-1-2) edge[bend left=50] node[below]{$\ra$}(m-1-1);
\end{tikzpicture}\end{center} with the following properties

\begin{enumerate}
\item For any injective $I\in\D$, there is some object $C\in \C$ such that $I\cong\la(C)$.
\item $\la$ sends injectives to injectives.
\end{enumerate}

Then the adjunction $\la\dashv\ra$ can be extended to a recollement

\begin{center}\begin{tikzpicture}
\matrix (m) [ampersand replacement= \&,matrix of math nodes, row sep=5em, 
column sep=5em,text height=1.5ex,text depth=0.25ex] 
{\ker(\ra)\&\C\&\D\\}; 
\path[->,thick, font=\scriptsize]
(m-1-1) edge node[above]{$i$}(m-1-2)
(m-1-2) edge[bend right=50] node[above]{$i_L$}(m-1-1)
(m-1-2) edge[bend left=50] node[below]{$i_R$}(m-1-1)
(m-1-3)edge[bend left=50] node[below]{$R_0\ra$}(m-1-2)
(m-1-3) edge[bend right=50] node[above]{$\la$}(m-1-2)
(m-1-2) edge node[above]{$\ra$}(m-1-3);
\end{tikzpicture}\end{center} 
In particular, $\ra$ is exact.
\end{thm}

\section{\label{section4}Derived Recollements and Applications of the Theorem}

\begin{defn}A recollement  $(L_0\ra, \la, \ra)$ will be called a \dbf{left derived recollement} and a recollement $(\la,\ra,R_0\la)$ will be called a \dbf{right derived recollement}.  Both will be referred to as \dbf{derived recollements}.\end{defn}
\subsection{Yoneda embedding}
The first application of Theorem \ref{recthm} will be to the adjoint pair \begin{center}\begin{tikzpicture}
\matrix (m) [ampersand replacement= \&,matrix of math nodes, row sep=3em, 
column sep=6em,text height=1.5ex,text depth=0.25ex] 
{\fp(\A^{op},\ab)\&\A\\}; 
\path[->,thick, font=\scriptsize]
(m-1-1) edge[bend left=50] node[auto]{$w$}(m-1-2)
(m-1-2) edge[bend left=50] node[below]{$\Y$}(m-1-1);
\end{tikzpicture}\end{center}
where $\A$ is an abelian category with enough projectives. 

We verify
\begin{enumerate}
\item $\Y$ is fully faithful.
\item Every projective $P\in\fp(\A^{op},\ab)$ is representable and thus $P\cong \Y(X)$ for some $X$.
\item If $P$ is projective in $\A$, then $\Y(P)=(\blank,P)$ is projective in $\fp(\A ^{op},\ab)$.  
\end{enumerate}

Applying Theorem \ref{recthm} gives the following.

\begin{cor}[The Defect Recollement]If $\A$ is abelian and has enough projectives then there is a recollement \begin{center}\begin{tikzpicture}
\matrix (m) [ampersand replacement= \&,matrix of math nodes, row sep=5em, 
column sep=6em,text height=1.5ex,text depth=0.25ex] 
{\ker(w)\&\fp(\A^{op},\ab)\&\A\\}; 
\path[->,thick, font=\scriptsize]
(m-1-1) edge node[above]{$i$}(m-1-2)
(m-1-2) edge[bend right=50] node[above]{$(\blank)^0$}(m-1-1)
(m-1-2) edge[bend left=50] node[below]{$(\blank)_0$}(m-1-1)
(m-1-3)edge[bend left=50] node[below]{$\Y$}(m-1-2)
(m-1-3) edge[bend right=50] node[above]{$L_0\Y$}(m-1-2)
(m-1-2) edge node[above]{$w$}(m-1-3);
\end{tikzpicture}\end{center} 
\end{cor}
Note that when $\A=\Mod(R)$ the defect of a functor $F\in\fp(\A^{op},\ab)$ is simply given by $w(F)=F(R)$. If $\A=\mod(\C^{op})$ for an additive category $\C$ with weak cokernels, then $w(F)(C)=F(\hom_\C(C,-))$ for any $F\in\fp(\A^{op},\ab)$ and any $C\in\C$.

Also, note that the kernel $\ker(w)$ is equivalent to the category of short exact sequences in $\A$, where the morphisms are homotopy classes of chain maps. Therefore, this category is abelian (whereas, by a characterisation of monomorphisms and epimorphisms, the category of short exact sequences and chain maps (rather than homotopy classes of chain maps) is not abelian unless $\A\simeq 0$, because there are morphisms which are monic and epic but not invertible).
%Now we look at the adjoint pair \begin{tikzpicture}
%\matrix (m) [ampersand replacement= \&,matrix of math nodes, row sep=3em, 
%column sep=6em,text height=1.5ex,text depth=0.25ex] 
%{\Mod(R)\&(\mod(R^{op}),\ab)\\}; 
%\path[->,thick, font=\scriptsize]
%(m-1-1) edge[bend left=50] node[auto]{$t$}(m-1-2)
%(m-1-2) edge[bend left=50] node[below]{$\ev_R$}(m-1-1);
%\end{tikzpicture}\end{center}
\subsection{A recollement for finitely accessible categories}
Let $\C$ be a finitely accessible category.  We have the embedding $\h:\C\to((\fp\C)^{op},\ab):C\mapsto(\blank,C)$.

Suppose $\C$ is cocomplete. Define $v:((\fp\C)^{op},\ab)\to\C$, to be the unique functor which preserves colimits such that, for any $C\in\C$, there is an isomorphism $v((\blank,C))\cong C$ which is natural in $C$. One can show that $v$ is the left adjoint of $\h$.

Now assume that $\C$ is abelian with enough projectives. The projectives in $((\fp\C)^{op},\ab)$ are those functors of the form $(\blank,P)$ where $P$ is a pure projective. Using this, it is easy to verify the assumptions of the main theorem, and hence we have the following.

\begin{cor}If $\C$ is finitely accessible, cocomplete, abelian and has enough projectives then there is a recollement \begin{center}\begin{tikzpicture}
\matrix (m) [ampersand replacement= \&,matrix of math nodes, row sep=5em, 
column sep=6em,text height=1.5ex,text depth=0.25ex] 
{\ker(v)\&((\fp\C)^{op},\ab)\&\C\\}; 
\path[->,thick, font=\scriptsize]
(m-1-1) edge node[above]{$i$}(m-1-2)
(m-1-2) edge[bend right=50] node[above]{$(\blank)^0$}(m-1-1)
(m-1-2) edge[bend left=50] node[below]{$(\blank)_0$}(m-1-1)
(m-1-3)edge[bend left=50] node[below]{$\h$}(m-1-2)
(m-1-3) edge[bend right=50] node[above]{$L_0\Y$}(m-1-2)
(m-1-2) edge node[above]{$v$}(m-1-3);
\end{tikzpicture}\end{center} 
\end{cor}
\begin{rem}We take this opportunity to point out that our main theorem in fact holds in greater generality. In the notation of the theorem, $\D$ need not be an abelian category, but it could just be an exact category with enough projectives, with the $0$-th derived functor read as being with respect to those projectives. (Note that the definition of recollement still makes sense in this broader setting.) It follows that, in this example, we needn't have assumed that our finitely accessible category $\C$ is abelian with enough projectives: It is enough that we have the pure exact structure on $\C$, which is always present and always has enough pure projectives (see \cite{cb}).
\end{rem}
\subsection{A recollement for modules as tensor products}
The injectives in $(\mod(R^{op}),\ab)$ are those functors of the form $M\otimes\blank$ where $M$ is pure injective. It follows that, by our main theorem, we have the following.
\begin{cor}There is a recollement \begin{center}\begin{tikzpicture}
\matrix (m) [ampersand replacement= \&,matrix of math nodes, row sep=5em, 
column sep=6em,text height=1.5ex,text depth=0.25ex] 
{\ker(\ev_R)\&(\mod(R^{op}),\ab)\&\Mod(R)\\}; 
\path[->,thick, font=\scriptsize]
(m-1-1) edge node[above]{$i$}(m-1-2)
(m-1-2) edge[bend right=50] node[above]{$(\blank)^0$}(m-1-1)
(m-1-2) edge[bend left=50] node[below]{$(\blank)_0$}(m-1-1)
(m-1-3)edge[bend left=50] node[below]{$R_0t$}(m-1-2)
(m-1-3) edge[bend right=50] node[above]{$t$}(m-1-2)
(m-1-2) edge node[above]{$\ev_R$}(m-1-3);
\end{tikzpicture}\end{center} 
\end{cor}
This shows that $\ev_R$, evaluation at the ring, has a fully faithful right adjoint, and so by Theorem \ref{serloc} we have the following.
\begin{cor}For any ring $R$, 
\begin{displaymath}
\frac{(\mod(R^{op}),\ab)}{\{F:\ev_R (F)=0\}}\simeq \Mod(R).
\end{displaymath}
\end{cor}
Since every functor $F:\mod(R^{op})\to\ab$ has an injective copresentation, and injectives in $(\mod(R^{op}),\ab)$ are of the form $M\otimes\blank$ for a pure injective $M$, the kernel $\ker(\ev_R)$ is equivalent to the category of short exact sequences of pure injectives, where morphisms are homotopy classes of chain maps. Therefore, this category is abelian.
\subsection{A recollement for modules over coherent rings.}Let $\C$ be an additive category with weak kernels (or a coherent ring). The injectives in $\fp(\mod(\C^{op}),\ab)$ are those functors of the form $M\otimes \blank$ for a finitely presented module $M\in\mod(\C)$. Applying our main theorem, we have the following.

\begin{cor}There is a recollement \begin{center}\begin{tikzpicture}
\matrix (m) [ampersand replacement= \&,matrix of math nodes, row sep=5em, 
column sep=6em,text height=1.5ex,text depth=0.25ex] 
{\ker(\ev_\C)\&\fp(\mod(\C^{op}),\ab)\&\mod(\C)\\}; 
\path[->,thick, font=\scriptsize]
(m-1-1) edge node[above]{$i$}(m-1-2)
(m-1-2) edge[bend right=50] node[above]{$(\blank)^0$}(m-1-1)
(m-1-2) edge[bend left=50] node[below]{$(\blank)_0$}(m-1-1)
(m-1-3)edge[bend left=50] node[below]{$R_0t$}(m-1-2)
(m-1-3) edge[bend right=50] node[above]{$t$}(m-1-2)
(m-1-2) edge node[above]{$\ev_\C$}(m-1-3);
\end{tikzpicture}\end{center} 
where $\ev_\C (F)(C)=F(\hom_C(C,-))$ for any $F\in(\mod(\C^{op}),\ab)$ and $C\in\C$.

The category $\fp(\mod(\C^{op}),\ab)$ has enough injectives. The injectives are functors of the form $M\otimes\blank$ for some $M\in\mod(\C)$. It follows that the kernel, $\ker(\ev_\C)$, is equivalent to the category of short exact sequences in $\mod(\C)$, where morphisms are homotopy classes of chain maps. Therefore, this category is abelian.

From Proposition \ref{serloc} we have the equivalence
\begin{displaymath}
\frac{\fp(\mod(\C^{op}),\ab)}{\{F:\ev_\C(F)=0\}}\simeq \mod(\C).
\end{displaymath}
We now apply this to find a result about dualising $k$-varieties which generalises \cite[Proposition 2.5]{dualvar}.
\begin{cor}Let $k$ be a commutative artinian ring. If $\C$ is a finite $k$-variety category with weak kernels and weak cokernels such that $\mod(\C^{op})$ is a dualising variety then $\C$ is also a dualising variety.
\end{cor}
\begin{pf}
Since $\C$ has weak kernels and weak cokernels, both $\mod(\C)$ and $\mod(\C^{op})$ are abelian.

Let $\A=\mod(\C^{op})$. The defect $w:\fp(\A^{op},\ab)\to \A$ is given by $w(F)(C)=F(\hom_\C(C,-))$ for all $F\in\fp(\A^{op},\ab)$ and $C\in\C$. By the definition of $\ev_\C$, it is clear that, under the standard equivalence $\fp(\A,\ab)\simeq(\fp(\A^{op},\ab))^{op}$, the Serre subcategory $\{F:\ev_\C(F)=0\}$ corresponds to $\{G:w(G)=0\}$.

By this and the equivalence $\fp(\A^{op},\ab)/\{G:w(G)=0\}\simeq\A$ established in \cite[Section 3]{coh}, we have
\begin{displaymath}
\mod(\C)\simeq\frac{\fp(\A,\ab)}{\{F:\ev_\C(F)=0\}}\simeq \frac{(\fp (\A^{op},\ab))^{op}}{\{G:w(G)=0\}}\simeq \mod(\C^{op})^{op}.
\end{displaymath}
One can easily show that this equivalence is isomorphic to the functor $M\mapsto \hom_k(M(\blank),k/J)$ where $J$ is the Jacobson radical of $k$. This shows that $\C$ is a dualising variety.
\end{pf}
\end{cor}
\subsection{A \textit{non}-example}
We will show that the tensor embedding $t:\mod(R)\to\fp(\Mod(R^{op}),\ab):M\mapsto M\otimes\blank$ does not admit a right adjoint. This means that the process carried out in previous examples cannot even begin for this embedding. Assume for contradiction that $t$ has a right adjoint $\ra$. Let $\{M_\lambda\}_{\lambda\in\Lambda}$ be an indexed family of finitely presented right $R$-modules. Then, since $\fp(\Mod(R^{op}),\ab)$ has products, and $\ra$ must preserve products, $\ra(\prod_{\lambda}t(M_\lambda))$ is the product of the $M_\lambda$ in $\mod(R)$. Therefore $\mod(R)$ has arbitrary products. In particular, $\mod(R)$ has arbitrary powers. We will now explain why this is a contradiction.

We recall an argument due to Freyd which shows that a small category $\C$ with arbitrary powers must have $|\C(A,B)|\leqslant 1$ for any $A,B\in\C$. Since $\C$ is small, each hom-set in $\C$ must have cardinality at most $|X|$, where $X$ is the set of all morphisms in $\C$. Suppose $\C$ has arbitrary powers and assume for contradiction that $|\C(A,B)|\geqslant 2$ for some $A,B\in\C$. Then $|\C(A,B^X)|=|\C(A,B)^X|=|\C(A,B)|^{|X|}\geqslant 2^{|X|}>|X|$, a contradiction. In particular, if $\C$ is a small pre-additive category with arbitrary powers then $\C\simeq 0$. 
\subsection{The Auslander-Gruson-Jenson Recollement}
In \cite{dr1} the following result was established. We can now view this as a corollary of our main theorem. The functor $\agj$ can be described by the explicit formula $\agj(F)(M)=\nat(F,\blank\otimes M)$ (one can verify this by using the fact that $\blank\otimes M$ is injective in $\fp(\Mod(R),\ab)$). Its right adjoint, $\agr$, can be described by the similar formula $\agr(G)(N)=\nat(G,N\otimes-)$.

The projectives in $\fp(\Mod(R),\ab)$ are precisely the representable functors, and the projectives in $(\mod(R^{op}),\ab)^{op}$ are the functors of the form $N\otimes \blank$ for a pure injective $N\in\Mod(R)$. One can use this information to verify the assumptions in our main theorem and establish the following.
\begin{cor}[Auslander-Gruson-Jensen Recollement] There is a recollement \begin{center}\begin{tikzpicture}
\matrix (m) [ampersand replacement= \&,matrix of math nodes, row sep=5em, 
column sep=3em,text height=1.5ex,text depth=0.25ex] 
{\ker(\agj)\&\fp(\Mod(R),\ab)\&(\mod(R^{op}),\ab)^{op}\\}; 
\path[->,thick, font=\scriptsize]
(m-1-1) edge node[above]{$i$}(m-1-2)
(m-1-2) edge[bend right=50] node[above]{$(\blank)^p$}(m-1-1)
(m-1-2) edge[bend left=50] node[below]{$(\blank)_q$}(m-1-1)
(m-1-3)edge[bend left=50] node[below]{$\agr$}(m-1-2)
(m-1-3) edge[bend right=50] node[above]{$\agl=L_0(\agr)$}(m-1-2)
(m-1-2) edge node[above]{$\agj$}(m-1-3);
\end{tikzpicture}\end{center} 
\end{cor}
The kernel, $\ker(\agj)$, is dual to the category of pure exact sequences in $\Mod(R)$, where morphisms are homotopy classes of chain maps. Therefore, this category is abelian.
\subsection{The Restriction Recollement}
In \cite{kraus}, Krause shows that any finitely accessible additive category has weak kernels. Let $\C$ be a finitely accessible additive category. Then $\fp(\C^{op},\ab)$ is abelian by Krause's result. 

Any functor $G\in((\fp\C)^{op},\ab)$ has can be extended to a unique functor $\overleftarrow{G}:\C^{op}\to\ab$ which preserves inverse limits. Krause also shows that for any functor $G\in((\fp\C)^{op},\ab)$, its extension $\overleftarrow{G}$ is finitely presented. Therefore there is a functor $\overleftarrow{-}:((\fp\C)^{op},\ab)\to\fp(\C^{op},\ab)$, and this is the right adjoint of the restriction functor $\res:\fp(\C^{op},\ab)\to((\fp\C)^{op},\ab)$. Note that Krause observed the existence of the left adjoint to $\res$ in \cite{kraus}, describing it as the unique functor $((\fp\C)^{op},\ab)\to\fp(\C^{op},\ab)$ which preserves colimits and sends a representable $\hom_{\fp \C}(-,X)$ to $\hom_\C(-,X)$. However, for this argument, one must take note that $\fp(\C^{op},\ab)$ actually has all colimits. Indeed it has cokernels, but the coproducts are unexpected, as the Yoneda embedding $\C\to\fp(\C^{op},\ab)$ preserves them. As an example of our main theorem, we give an alternative method for constructing a left adjoint.

Any projective in $\fp(\C^{op},\ab)$ is of the form $\hom_\C(\blank,C)$ for some $C\in\C$. This is clearly isomorphic to $\overleftarrow{(\blank,C)}$. The projectives in $((\fp\C)^{op},\ab)$ are those of the form $(\blank,P)$ where $P$ is a pure projective, which extends to $(\blank,P)$. This verifies the assumptions of our main theorem, and we have the following.
\begin{cor}[The Restriction Recollement]There is a recollement \begin{center}\begin{tikzpicture}
\matrix (m) [ampersand replacement= \&,matrix of math nodes, row sep=5em, 
column sep=3em,text height=1.5ex,text depth=0.25ex] 
{\ker(\res)\&\fp(\C^{op},\ab)\&((\fp\C)^{op},\ab)^{op}\\}; 
\path[->,thick, font=\scriptsize]
(m-1-1) edge node[above]{$i$}(m-1-2)
(m-1-2) edge[bend right=50] node[above]{$(\blank)^0$}(m-1-1)
(m-1-2) edge[bend left=50] node[below]{$(\blank)_0$}(m-1-1)
(m-1-3)edge[bend left=50] node[below]{$\overleftarrow{-}$}(m-1-2)
(m-1-3) edge[bend right=50] node[above]{$L_0(\overleftarrow{-})$}(m-1-2)
(m-1-2) edge node[above]{$\res$}(m-1-3);
\end{tikzpicture}\end{center} 
\end{cor}

If $\C$ also has kernels (see \cite{cb} for necessary and sufficient conditions for this), then $\ker(\res)$ is  equivalent to the category of pure exact sequences, where morphisms are homotopy classes of chain maps. Therefore, this category is abelian.

\section{\label{section5}Defect Zero Functors and Duality}

Around the same time that Auslander published his landmark paper on coherent functors, which many now call finitely presented functors, his students Janet Fisher-Palmquist and David Newell were researching how to extend many of the familiar module theoretic concepts to appropriate functor theoretic concepts.  The idea was that if one can view a functor as a generalization of a module, then any property that the module category $\Mod(R)$ exhibits may have a functor analog.

In particular Fisher-Palmquist and Newell, like Auslander, were very interested in bringing methods from homological algebra to functor categories.  This is evident in the following series of papers.
\begin{enumerate}
\item In \cite{fishertensor} Fisher discusses properties of a tensor product of functors 
$F\otimes G$ which is a generalization of the tensor product of modules.  
\item  In \cite{fisheradjoint} Fisher-Palmquist and Newell establish that the analog of a bi-module at the level of functors is the notion of a bifunctor.  The analog of the bi-module $R$ is the bifunctor $\hom$.  As such the dual $\hom(\blank, R)$ has an analog for functors which is the functor $(\blank)^*$ described below.   One of the main results from \cite{fisheradjoint} is that adjoints between  functor categories arise as the analog of the hom-tensor adjunction.  
\item In \cite{newell}, Newell extends the well known Morita Theorems for module categories to functor categories.
\end{enumerate}

The following is the construction of Fisher-Palmquist and Newell appearing in \cite{fisheradjoint}.  Given a covariant functor $F\:\A\to \ab$, there is a contravariant functor $F^*\:\A\To \ab$ defined by $$F^*(X)=\nat\big(F,(X,\blank)\big)$$ Similarly given a contravariant functor $G\:\A\To \ab$ there is a covariant functor $G^*\:\A\To \ab$ defined by $$G^*(X)=\nat\big(G,(\blank, X)\big)$$   As stated above, the functor $(\blank)^*$ can be viewed as the functor analog of $\hom(\blank, R)$ and in fact the notation $\nat(\blank, \hom)$ could be used in place of $(\blank)^*$ since the bi-functor $\hom$ takes on the role of the bi-module $R$.    An easily verified fact is that for any representable functor $(X,\blank)$, the dual $(X,\blank)^*\cong (\blank, X)$ and $(\blank, X)^*\cong (X,\blank)$.

In \cite{defect}, the dual $(\blank)^*$ is studied at the level of the functor categories $\fp(\A,\ab)$ and $\fp(\A^{op},\ab)$.   The following \fp-dual formula is established there. For any finitely presented functor $F\:\A\to \ab$  $$F^*\cong\big(\blank, w(F)\big).$$  Similarly for any contravariant finitely presented functor $G\:\A\to \ab$ $$G^*\cong \big(v(G),\blank\big)$$  Therefore, the defect completely determines the dual of a finitely presented functor, the dual of any finitely presented functor is representable,     and the assignment $(\blank)^*$ gives rise to a pair of left exact contravariant functors
\begin{center}\begin{tikzpicture}
\matrix (m) [ampersand replacement= \&,matrix of math nodes, row sep=3em, 
column sep=5em,text height=1.5ex,text depth=0.25ex] 
{\fp(\A,\ab)\&\fp(\A^{op},\ab)\\}; 
\path[->,thick, font=\scriptsize]
(m-1-1) edge[bend left=50] node[auto]{$(\blank)^*$}(m-1-2)
(m-1-2) edge[bend left=50] node[below]{$(\blank)^*$}(m-1-1);
\end{tikzpicture}\end{center}

If $\A$ has both enough injectives and projectives, then by Gentle's result, i.e. Proposition \ref{gentleinj}, the categories $\fp(\A,\ab)$ and $\fp(\A^{op},\ab)$ have enough injectives.    In this case, one can compute the left derived functors $L^k(\blank)^*$.  This is done in \cite{defect}, where it is shown that $n=L^0(\blank)^*$ is the only non-trivial left derived functor and  \begin{center}\begin{tikzpicture}
\matrix (m) [ampersand replacement= \&,matrix of math nodes, row sep=3em, 
column sep=5em,text height=1.5ex,text depth=0.25ex] 
{\fp(\A,\ab)\&\fp(\A^{op},\ab)\\}; 
\path[->,thick, font=\scriptsize]
(m-1-1) edge[bend left=50] node[auto]{$n=L^0(\blank)^*$}(m-1-2)
(m-1-2) edge[bend left=50] node[below]{$n=L^0(\blank)^*$}(m-1-1);
\end{tikzpicture}\end{center} is a pair of exact functors.\footnote{Technically, the argument given in \cite{defect} applies to $\A=\Mod(R)$ but works in complete generality for any abelian category $\A$ which has both enough projectives and injectives.}

On the other hand for \textit{any} abelian category $\A$ the functor category $\fp(\A,\ab)$ always has enough projectives.  Hence the right derived functors \begin{center}\begin{tikzpicture}
\matrix (m) [ampersand replacement= \&,matrix of math nodes, row sep=3em, 
column sep=5em,text height=1.5ex,text depth=0.25ex] 
{\fp(\A,\ab)\&\fp(\A^{op},\ab)\\}; 
\path[->,thick, font=\scriptsize]
(m-1-1) edge[bend left=50] node[auto]{$R_k(\blank)^*$}(m-1-2)
(m-1-2) edge[bend left=50] node[below]{$R_k(\blank)^*$}(m-1-1);
\end{tikzpicture}\end{center} can always be computed for an arbitrary abelian category $\A$.    

Let $F$ be any finitely presented functor and take any projective resolution.

\begin{center}\begin{tikzpicture}
\matrix(m)[ampersand replacement=\&, matrix of math nodes, row sep=3em, column sep=2.5em, text height=1.5ex, text depth=0.25ex]
{0\&(Z,\blank)\&(Y,\blank)\&(X,\blank)\&F\&0\\}; 
\path[->,thick, font=\scriptsize]
(m-1-1) edge node[auto]{}(m-1-2)
(m-1-2) edge node[auto]{$(g,\blank)$}(m-1-3)
(m-1-3) edge node[auto]{$(f,\blank)$}(m-1-4)
(m-1-4) edge node[auto]{}(m-1-5)
(m-1-5) edge node[auto]{}(m-1-6);
\end{tikzpicture}\end{center}  Applying the functor $(\blank)^*$ to the projective resolution and using the fact that for representables  $(A,\blank)^*\cong(\blank, A)$ yields the following complex.

\begin{center}\begin{tikzpicture}[framed]
\matrix(m)[ampersand replacement=\&, matrix of math nodes, row sep=3em, column sep=2.5em, text height=1.5ex, text depth=0.25ex]
{0\&(\blank,X)\&(\blank,Y)\&(\blank,Z)\&0\\}; 
\path[->,thick, font=\scriptsize]
(m-1-1) edge node[auto]{}(m-1-2)
(m-1-2) edge node[auto]{$(\blank,f)$}(m-1-3)
(m-1-3) edge node[auto]{$(\blank,g)$}(m-1-4)
(m-1-4) edge node[auto]{}(m-1-5);
\end{tikzpicture}\end{center}  From here on we will use the notation $$\extd_k:=R_k(\blank)^*.$$  By definition, the homology of this complex is $\extd_k(F)$.  We mention that $\extd_0\cong (\blank)^*$ which follows from the left exactness of $(\blank)^*$.  Also, if $k\ge 3$, then $\extd_k\cong 0$.  This leaves $\extd_1$ and $\extd_2$.  

\subsection*{Computation of $R_1(\blank)^*$} To compute $\extd_1(F)$ we must compute the homology of the following complex at the position marked by $\bullet$.  
\begin{center}\begin{tikzpicture}[framed]
\matrix(m)[ampersand replacement=\&, matrix of math nodes, row sep=3em, column sep=2.5em, text height=1.5ex, text depth=0.25ex]
{0\&(\blank,X)\&\underset{\bullet}{(\blank,Y)}\&(\blank,Z)\&0\\}; 
\path[->,thick, font=\scriptsize]
(m-1-1) edge node[auto]{}(m-1-2)
(m-1-2) edge node[auto]{$(\blank,f)$}(m-1-3)
(m-1-3) edge node[auto]{$(\blank,g)$}(m-1-4)
(m-1-4) edge node[auto]{}(m-1-5);
\end{tikzpicture}\end{center}  By definition $$\extd_1(F)\cong\frac{\ker(\blank, g)}{\im(\blank, f)}.$$  This can computed by as follows.  Embedding the following commutative diagram with exact rows into the category $\fp(\A^{op},\ab)$

\begin{center}\begin{tikzpicture}[framed]
\matrix(m)[ampersand replacement=\&, matrix of math nodes, row sep=1.5em, column sep=1.5em, text height=1.5ex, text depth=0.25ex]
{\&\&0\\
\&w(F)\&w(F)\\
0\&w(F)\&X\&Y\&Z\&0\\
\&0\&V\&Y\&Z\&0\\
\&\&0\\}; 
\path[->,thick, font=\scriptsize]
(m-1-3)edge (m-2-3)
(m-2-2) edge node[left]{$1$}(m-3-2)
(m-2-3) edge node[left]{$k$}(m-3-3)
(m-2-2) edge node[above]{$1$}(m-2-3)
(m-3-1) edge node[auto]{}(m-3-2)
(m-3-2) edge node[auto]{$k$}(m-3-3)
(m-3-3) edge node[auto]{$f$}(m-3-4)
(m-3-4) edge node[auto]{$g$}(m-3-5)
(m-3-5) edge node[auto]{}(m-3-6)
(m-4-2) edge node[auto]{}(m-4-3)
(m-4-3) edge node[auto]{$m$}(m-4-4)
(m-4-4) edge node[auto]{$g$}(m-4-5)
(m-4-5) edge node[auto]{}(m-4-6)
(m-3-3) edge node[left]{$e$}(m-4-3)
(m-3-4) edge node[left]{$1$}(m-4-4)
(m-3-5) edge node[left]{$1$}(m-4-5)
(m-4-3) edge (m-5-3);
\end{tikzpicture}\end{center}

yields the following commutative diagram with exact rows and columns

\begin{center}\begin{tikzpicture}[framed]
\matrix(m)[ampersand replacement=\&, matrix of math nodes, row sep=2em, column sep=2em, text height=1.5ex, text depth=0.25ex]
{\&\&0\\
\&\big(\blank,w(F)\big)\&\big(\blank, w(F)\big)\\
0\&\big(\blank, w(F)\big)\&(\blank, X)\&(\blank, Y)\&N\&0\\
\&0\&(\blank, V)\&(\blank, Y)\&(\blank, Z)\\
\&\&\extd_1(F)\\
\&\&0\\}; 
\path[->,thick, font=\scriptsize]
(m-1-3)edge (m-2-3)
(m-2-2) edge node[left]{$1$}(m-3-2)
(m-2-3) edge node[left]{$(\blank, k)$}(m-3-3)
(m-2-2) edge node[above]{$1$}(m-2-3)
(m-3-1) edge node[auto]{}(m-3-2)
(m-3-2) edge node[auto]{$(\blank, k)$}(m-3-3)
(m-3-3) edge node[auto]{$(\blank, f)$}(m-3-4)
(m-3-4) edge node[auto]{$\kappa$}(m-3-5)
(m-3-5) edge node[auto]{}(m-3-6)
(m-4-2) edge node[auto]{}(m-4-3)
(m-4-3) edge node[auto]{$(\blank, m)$}(m-4-4)
(m-4-4) edge node[auto]{$(\blank, g)$}(m-4-5)
%(m-4-5) edge node[auto]{}(m-4-6)
(m-3-3) edge node[left]{$(\blank, e)$}(m-4-3)
(m-3-4) edge node[left]{$1$}(m-4-4)
(m-3-5) edge node[left]{$\zeta$}(m-4-5)
(m-4-3) edge (m-5-3)
(m-5-3) edge (m-6-3);
\end{tikzpicture}\end{center}  Therefore for any functor $F$ with presentation as given above, there is a presentation \begin{center}\begin{tikzpicture}[framed]
\matrix(m)[ampersand replacement=\&, matrix of math nodes, row sep=3em, column sep=2.5em, text height=1.5ex, text depth=0.25ex]
{0\&\big(\blank,w(F)\big)\&(\blank,X)\&(\blank,V)\&\extd_1(F)\&0\\}; 
\path[->,thick, font=\scriptsize]
(m-1-1) edge node[auto]{}(m-1-2)
(m-1-2) edge node[auto]{}(m-1-3)
(m-1-3) edge node[auto]{$(\blank, e)$}(m-1-4)
(m-1-4) edge node[auto]{}(m-1-5)
(m-1-5) edge node[auto]{}(m-1-6);
\end{tikzpicture}\end{center} In particular, it is easily seen that if $w(F)=0$, then $\extd_1(F)\cong 0$ as in this case $e$ is an isomorphism.  
\subsection*{Computation of $R_2(\blank)^*$}  To compute $\extd_2(F)$ we must compute the homology of the following complex at the position marked by $\bullet$.  
\begin{center}\begin{tikzpicture}[framed]
\matrix(m)[ampersand replacement=\&, matrix of math nodes, row sep=3em, column sep=2.5em, text height=1.5ex, text depth=0.25ex]
{0\&(\blank,X)\&(\blank,Y)\&\underset{\bullet}{(\blank,Z)}\&0\\}; 
\path[->,thick, font=\scriptsize]
(m-1-1) edge node[auto]{}(m-1-2)
(m-1-2) edge node[auto]{$(\blank,f)$}(m-1-3)
(m-1-3) edge node[auto]{$(\blank,g)$}(m-1-4)
(m-1-4) edge node[auto]{}(m-1-5);
\end{tikzpicture}\end{center} By definition $$\extd_2(F)\cong\frac{(\blank, Z)}{\im(\blank, g)}.$$  A straightforward homological argument shows that this is precisely the cokernel of the morphism $(\blank, g)$.

\begin{prop}The composition $\extd_2\circ \extd_2\cong(\blank)_0$.\end{prop}

\begin{pf}From the projective resolution of $F\in\fp(\A,\ab)$  \begin{center}\begin{tikzpicture}
\matrix(m)[ampersand replacement=\&, matrix of math nodes, row sep=1em, column sep=2.5em, text height=1.5ex, text depth=0.25ex]
{0\&(Z,\blank)\&(Y,\blank)\&(X,\blank)\&F\&0\\}; 
\path[->,thick, font=\scriptsize]
(m-1-1) edge node[auto]{}(m-1-2)
(m-1-2) edge node[auto]{$(g,\blank)$}(m-1-3)
(m-1-3) edge node[auto]{$(f,\blank)$}(m-1-4)
(m-1-4) edge node[auto]{$\alpha$}(m-1-5)
(m-1-5) edge node[auto]{}(m-1-6);
\end{tikzpicture}\end{center} one has exact the following commutative diagram with exact rows and columns

\begin{center}\begin{tikzpicture}
\matrix(m)[ampersand replacement=\&, matrix of math nodes, row sep=1.5em, column sep=1.5em, text height=1.5ex, text depth=0.25ex]
{\&\&0\\
\&w(F)\&w(F)\\
0\&w(F)\&X\&Y\&Z\&0\\
\&0\&V\&Y\&Z\&0\\
\&\&0\\}; 
\path[->,thick, font=\scriptsize]
(m-1-3)edge (m-2-3)
(m-2-2) edge node[left]{$1$}(m-3-2)
(m-2-3) edge node[left]{$k$}(m-3-3)
(m-2-2) edge node[above]{$1$}(m-2-3)
(m-3-1) edge node[auto]{}(m-3-2)
(m-3-2) edge node[auto]{$k$}(m-3-3)
(m-3-3) edge node[auto]{$f$}(m-3-4)
(m-3-4) edge node[auto]{$g$}(m-3-5)
(m-3-5) edge node[auto]{}(m-3-6)
(m-4-2) edge node[auto]{}(m-4-3)
(m-4-3) edge node[auto]{$m$}(m-4-4)
(m-4-4) edge node[auto]{$g$}(m-4-5)
(m-4-5) edge node[auto]{}(m-4-6)
(m-3-3) edge node[left]{$e$}(m-4-3)
(m-3-4) edge node[left]{$1$}(m-4-4)
(m-3-5) edge node[left]{$1$}(m-4-5)
(m-4-3) edge (m-5-3);
\end{tikzpicture}\end{center}  Therefore by embedding the short exact sequence

\begin{center}\begin{tikzpicture}
\matrix(m)[ampersand replacement=\&, matrix of math nodes, row sep=3em, column sep=2.5em, text height=1.5ex, text depth=0.25ex]
{0\&V\&Y\&Z\&0\\}; 
\path[->,thick, font=\scriptsize]
(m-1-1) edge (m-1-2)
(m-1-2) edge node[above]{$m$}(m-1-3)
(m-1-3) edge node[auto]{$g$}(m-1-4)
(m-1-4) edge node[auto]{}(m-1-5);
\end{tikzpicture}\end{center}   into $\fp(\A^{op},\ab)$ and using the fact that $\coker((\blank,g))\cong \extd_2(F)$ we have the following exact sequence

\begin{center}\begin{tikzpicture}[framed]
\matrix(m)[ampersand replacement=\&, matrix of math nodes, row sep=3em, column sep=2.5em, text height=1.5ex, text depth=0.25ex]
{0\&(\blank, V)\&(\blank,Y)\&(\blank,Z)\&\extd_2(F)\&0\\}; 
\path[->,thick, font=\scriptsize]
(m-1-1) edge (m-1-2)
(m-1-2) edge node[above]{$(\blank, m)$}(m-1-3)
(m-1-3) edge node[auto]{$(\blank,g)$}(m-1-4)
(m-1-4) edge node[auto]{}(m-1-5)
(m-1-5) edge node[auto]{}(m-1-6);
\end{tikzpicture}\end{center}  which is a projective presentation of $\extd_2(F)$.  Now using this projective presentation we have that $\extd_2(\extd_2(F))\cong \coker((m,\blank))$ thereby yielding the following projective presentation

\begin{center}\begin{tikzpicture}
\matrix(m)[ampersand replacement=\&, matrix of math nodes, row sep=3em, column sep=2.5em, text height=1.5ex, text depth=0.25ex]
{(Y,\blank)\&(V,\blank)\&\extd_2(\extd_2(F))\&0\\}; 
\path[->,thick, font=\scriptsize]
(m-1-1) edge node[auto]{$(m,\blank)$}(m-1-2)
(m-1-2) edge node[auto]{}(m-1-3)
(m-1-3) edge node[auto]{}(m-1-4);
\end{tikzpicture}\end{center} But from the construction of the defect sequence in Section \ref{section2} one sees that $$\extd_2(\extd_2(F))\cong F_0$$ and it is easily verified that this isomorphism is natural in $F$ establishing that $$\extd_2\circ \extd_2\cong (\blank)_0$$   This completes the proof. $\qed$
\end{pf}

Notice that in this proof we also established that for any finitely presented functor with projective resolution \begin{center}\begin{tikzpicture}
\matrix(m)[ampersand replacement=\&, matrix of math nodes, row sep=1em, column sep=2.5em, text height=1.5ex, text depth=0.25ex]
{0\&(Z,\blank)\&(Y,\blank)\&(X,\blank)\&F\&0\\}; 
\path[->,thick, font=\scriptsize]
(m-1-1) edge node[auto]{}(m-1-2)
(m-1-2) edge node[auto]{$(g,\blank)$}(m-1-3)
(m-1-3) edge node[auto]{$(f,\blank)$}(m-1-4)
(m-1-4) edge node[auto]{$\alpha$}(m-1-5)
(m-1-5) edge node[auto]{}(m-1-6);
\end{tikzpicture}\end{center} there is a projective resolution 

\begin{center}\begin{tikzpicture}[framed]
\matrix(m)[ampersand replacement=\&, matrix of math nodes, row sep=3em, column sep=2.5em, text height=1.5ex, text depth=0.25ex]
{0\&(\blank,V)\&(\blank,Y)\&(\blank,Z)\&\extd_2(F)\&0\\}; 
\path[->,thick, font=\scriptsize]
(m-1-1)edge(m-1-2)
(m-1-2) edge node[above]{$(\blank,m)$}(m-1-3)
(m-1-3) edge node[auto]{$(\blank,g)$}(m-1-4)
(m-1-4) edge node[auto]{}(m-1-5)
(m-1-5) edge node[auto]{}(m-1-6);
\end{tikzpicture}\end{center} Moreover, if $w(F)\cong0$ then $V=X$ and both $F$ and $\extd_2(F)$ arise from the same short exact sequence.

By definition $\extd_k(F)$ is the functor that can be computed as follows.  Given any $X\in \A$, $\extd_k(F)$ evaluated at $X$ will be $$\ext^k\big(F,(X,\blank)\big)$$   Suppose that we have a short exact sequence of functors in $\fp(\A,\ab)$
\begin{center}\begin{tikzpicture}
\matrix(m)[ampersand replacement=\&, matrix of math nodes, row sep=3em, column sep=2.5em, text height=1.5ex, text depth=0.25ex]
{0\&F\&G\&H\&0\\}; 
\path[->,thick, font=\scriptsize]
(m-1-1) edge node[auto]{}(m-1-2)
(m-1-2) edge node[auto]{}(m-1-3)
(m-1-3) edge node[auto]{}(m-1-4)
(m-1-4) edge node[auto]{}(m-1-5);
\end{tikzpicture}\end{center}  The long homology exact sequence yields the following exact sequence.

\begin{center}\begin{tikzpicture}
\matrix(m)[ampersand replacement=\&, matrix of math nodes, row sep=3em, column sep=2.5em, text height=1.5ex, text depth=0.25ex]
{\extd_1(F)\&\extd_2(H)\&\extd_2(G)\&\extd_2(F)\&\extd_3(H)\\}; 
\path[->,thick, font=\scriptsize]
(m-1-1) edge node[auto]{}(m-1-2)
(m-1-2) edge node[auto]{}(m-1-3)
(m-1-3) edge node[auto]{}(m-1-4)
(m-1-4) edge node[auto]{}(m-1-5);
\end{tikzpicture}\end{center} Since $\extd_3(F)\cong 0$ it follows that we have the following exact sequence \begin{center}\begin{tikzpicture}
\matrix(m)[ampersand replacement=\&, matrix of math nodes, row sep=3em, column sep=2.5em, text height=1.5ex, text depth=0.25ex]
{\extd_1(F)\&\extd_2(H)\&\extd_2(G)\&\extd_2(F)\&0\\}; 
\path[->,thick, font=\scriptsize]
(m-1-1) edge node[auto]{}(m-1-2)
(m-1-2) edge node[auto]{}(m-1-3)
(m-1-3) edge node[auto]{}(m-1-4)
(m-1-4) edge node[auto]{}(m-1-5);
\end{tikzpicture}\end{center} thereby establishing that $\extd_2$ is right exact.   As commented above if $w(F)\cong 0$, then $\extd_1(F)\cong 0$.   Hence on the Serre subcategory of defect zero functors $\extd_2$ is exact.   Since $\extd_2\circ \extd_2\cong (\blank)_0$, we have actually shown the following.

%\begin{prop}[Auslander,\cite{coh} Proposition 3.2]  For any finitely presented functor $F$, $w(F)\cong 0$, that is $F\in \fp_0(\Mod(R),\ab)$, if and only if  both
%\begin{enumerate}
%\item  $\ext^1(F,(X,\blank))\cong 0$.
%\item $\ext^0\big(F,(X,\blank)\big)\cong 0$.
%\end{enumerate}
%\end{prop}  

\begin{thm}\label{defectzerodual}For any abelian category $\A$, the functors $\extd_2=R_2(\blank)^*$ form a pair of right exact functors \begin{center}\begin{tikzpicture}
\matrix (m) [ampersand replacement= \&,matrix of math nodes, row sep=5em, 
column sep=4em,text height=1.5ex,text depth=0.25ex] 
{\fp(\A,\ab)\&\fp(\A^{op},\ab)\\}; 
\path[->,thick, font=\scriptsize]
(m-1-1)edge[bend right=50] node[below]{$\extd_2$}(m-1-2)
(m-1-2) edge[bend right=50] node[above]{$\extd_2$}(m-1-1);
%(m-1-2) edge[bend left] node[below]{$\agr$}(m-1-1);
\end{tikzpicture}\end{center}  satisfying the following properties
\begin{enumerate}

\item If $F$ has defect zero, then $\extd_2(F)$ has defect zero.
\item If $F$ arises from a short exact sequence, then $\extd_2(F)$ arises from the same short exact sequence.

\item $\extd_2$ restricts to a duality $\extd$

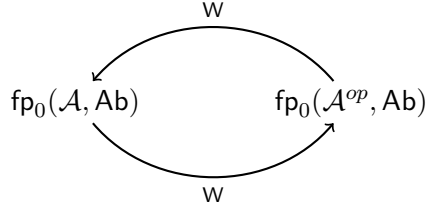
\begin{center}\begin{tikzpicture}
\matrix (m) [ampersand replacement= \&,matrix of math nodes, row sep=5em, 
column sep=4em,text height=1.5ex,text depth=0.25ex] 
{\fp_0(\A,\ab)\&\fp_0(\A^{op},\ab)\\}; 
\path[->,thick, font=\scriptsize]
(m-1-1)edge[bend right=50] node[below]{$\extd$}(m-1-2)
(m-1-2) edge[bend right=50] node[above]{$\extd$}(m-1-1);
%(m-1-2) edge[bend left] node[below]{$\agr$}(m-1-1);
\end{tikzpicture}\end{center}

\end{enumerate}
\end{thm}

The duality $\extd$ on these Serre subcategories $\fp_0(\A,\ab)$ and $\fp_0(\A^{op},\ab)$ implies that every short exact sequence induces two finitely presented functors which can be identified via $\extd$.  Both $F$ and $\extd_2(F)$ arise from the same short exact sequence.  This in some way explains why certain notions in representation theory concerning short exact sequences are symmetrical.  For example, if $\Lambda$ is an artin algebra, then the simple functors $S\:\mod(\Lambda)\to \ab$ are finitely presented as shown by Auslander in \cite{fart}.  Moreover, $w(S)=0$ and it can be shown that $S$ arises from an almost split sequence.  The functor $\extd(S)$ will be co-simple which turns out to be equivalent to being simple. The functor $\extd(S)$ arises from the same almost split sequence.   Similarly, if a finitely presented functor $F$ arises from a pure exact sequence, then its dual $\extd(F)$ arises from the same pure exact sequence.  We now state the promised formulas.

\begin{thm}\label{argen}[The Generalised Auslander-Reiten Formulas] Assume that $\A$ is abelian.
\begin{enumerate}
\item If $\A$ has enough projectives then 
\begin{enumerate} 
\item $\fp_0(\A,\ab)$ has enough injectives of the form  $\ext^1(A,\blank)$.
\item $\fp_0(\A^{op},\ab)$ has enough projectives of the form $\underline{\hom}(\blank, A)$.  
\item For any $C\in \A$ there are natural isomorphisms  $$\underline{\hom}(\blank, C)\cong \extd \ext^1(C,\blank)$$
\end{enumerate}
\item If $\A$ has enough injectives then
\begin{enumerate}
\item $\fp_0(\A^{op},\ab)$ has enough injectives of the form $\ext^1(\blank, A)$.
\item $\fp_0(\A,\ab)$ has enough projectives of the form $\overline{\hom}(A,\blank)$.
\item For any $C\in \A$ there is a natural isomorphism $$\overline{\hom}(C,\blank)\cong \extd \ext^1(\blank, C)$$
\end{enumerate}
\end{enumerate}
%In $\fp_0(\A,\ab)$ the projectives are precisely the functors $\underline{\Y}(X)\cong\overline{\hom}(\blank, X)$ appearing in the adjunction sequence
%
%
%
%\begin{center}\begin{tikzpicture}
%\matrix(m)[ampersand replacement=\&, matrix of math nodes, row sep=3em, column sep=2.5em, text height=1.5ex, text depth=0.25ex]
%{L^0\Y(X)\&(X,\blank)\&\underline{\Y}(X)\&0\\}; 
%\path[->,thick, font=\scriptsize]
%(m-1-1) edge node[auto]{}(m-1-2)
%(m-1-2) edge node[auto]{}(m-1-3)
%(m-1-3) edge node[auto]{}(m-1-4);
%\end{tikzpicture}\end{center}
%
%In $\fp_0(\A^{op},\ab)$ they are precisely the functors $\underline{\Y}(X)\cong \underline{\hom}(\blank, X)$ appearing in the adjunction sequence
%
%\begin{center}\begin{tikzpicture}
%\matrix(m)[ampersand replacement=\&, matrix of math nodes, row sep=3em, column sep=2.5em, text height=1.5ex, text depth=0.25ex]
%{L_0\Y(X)\&(\blank,X)\&\underline{\Y}(X)\&0\\}; 
%\path[->,thick, font=\scriptsize]
%(m-1-1) edge node[auto]{}(m-1-2)
%(m-1-2) edge node[auto]{}(m-1-3)
%(m-1-3) edge node[auto]{}(m-1-4);
%\end{tikzpicture}\end{center}
%

\end{thm}

\begin{pf}Since (1) and (2) are completely dual statements we need only show (1).  Auslander shows in \cite{coh} that $\fp_0(\A,\ab)$ has enough injectives of the form $\ext^1(C,\blank)$ under the assumptions of (1).  The functor $\extd$ is a duality and hence sends injectives to projectives.  All that remains is to show that the duality $\extd$ sends the functor $\ext^1(C,\blank)$ to the functor $\underline{\hom}(\blank, C)$.
From the following syzygy sequence of $C$

\begin{center}\begin{tikzpicture}
\matrix(m)[ampersand replacement=\&, matrix of math nodes, row sep=3em, column sep=2.5em, text height=1.5ex, text depth=0.25ex]
{0\&\Omega C \&P\&C\&0\\}; 
\path[->,thick, font=\scriptsize]
(m-1-1) edge node[auto]{}(m-1-2)
(m-1-2) edge node[auto]{$f$}(m-1-3)
(m-1-3) edge node[auto]{$g$}(m-1-4)
(m-1-4) edge node[auto]{}(m-1-5);
\end{tikzpicture}\end{center}

we have the following projective presentation of $\ext^1(C,\blank)$

\begin{center}\begin{tikzpicture}
\matrix(m)[ampersand replacement=\&, matrix of math nodes, row sep=3em, column sep=2.5em, text height=1.5ex, text depth=0.25ex]
{0\&(C,\blank) \&(P,\blank)\&(\Omega C,\blank)\&\ext^1(C,\blank)\&0\\}; 
\path[->,thick, font=\scriptsize]
(m-1-1) edge node[auto]{}(m-1-2)
(m-1-2) edge node[auto]{$(g,\blank)$}(m-1-3)
(m-1-3) edge node[auto]{$(f,\blank)$}(m-1-4)
(m-1-4) edge node[auto]{}(m-1-5)
(m-1-5) edge node[auto]{}(m-1-6);
\end{tikzpicture}\end{center}  From our above discussion about $\extd_2$ we have that $$\extd(\ext^1(C,\blank))\cong \coker(\blank, g)$$ and so we have an exact sequence

\begin{center}\begin{tikzpicture}
\matrix(m)[ampersand replacement=\&, matrix of math nodes, row sep=3em, column sep=2.5em, text height=1.5ex, text depth=0.25ex]
{(\blank, P)\&(\blank, C) \&\extd\ext^1(C,\blank)\&0\\}; 
\path[->,thick, font=\scriptsize]
(m-1-1) edge node[auto]{$(\blank,g)$}(m-1-2)
(m-1-2) edge node[auto]{}(m-1-3)
(m-1-3) edge node[auto]{}(m-1-4);
\end{tikzpicture}\end{center}  Since $g$ is an epimorphism from the projective $P$ to $C$ one can easily show that a morphism $h\:X\To C$ is in the subgroup $\P(C,X)$ of morphisms factoring through some projective if and only if it is in the image of $(X,f)$.  Since $\extd \ext^1(C,\blank)(X)\cong \coker((X,f))$, $$\extd \ext^1(C,\blank)(X)\cong \hom(X,C)/\P(X,C)=\underline{\hom}(\blank,C)(X)$$  This isomorphism is also easily seen to be natural so in fact we have $$\underline{\hom}(\blank, C)\cong \extd \ext^1(C,\blank).\qquad\qed$$
\end{pf}

\begin{rem}Auslander-Reiten formulas have been studied not only for modules, but objects in any Grothendieck abelian category, by Krause. See \cite{kar} for this. Here our focus is on the duality $\extd$ which is obtained by computing derived functors of the dual $(\blank)^*$.  This functor will be used to relate the Auslander-Gruson-Jensen Recollement and the Restriction Recollement.   We leave it to the interested reader to investigate what the Generalised Auslander-Reiten formulas say for a particular choice of the abelian category $\A$.  At the end of the paper, we show that these formulas reduce to the original Auslander-Reiten formulas.  The fact that we are able to recover these formulas by looking at derived functors of $(\blank)^*$ is quite surprising.
\end{rem}
%In the case that $\A$ is abelian and has enough projectives, the category $\fp_0(\A,\ab)$ is equivalent to the category $\fp(\overline{\A},\ab)$ using the following fact.   

\begin{cor}Suppose that $\A$ is an abelian category.
\begin{enumerate}
\item If $\A$ has enough injectives, then $\fp_0(\A,\ab)\cong \fp(\overline{\A},\ab)$. 
\item It $\A$ has enough projectives, then $\fp_0(\A^{op},\ab)\cong \fp(\underline{\A}^{op},\ab)$. 
\item If $\A$ has enough injectives and projectives, then $\fp(\underline{\A}^{op},\ab)\cong \fp(\overline{\A},\ab)$.
\end{enumerate}
\end{cor}

\begin{pf}In this case by Theorem \ref{argen}, every functor $F\in \fp_0(\A,\ab)$ has a projective presentation 
\begin{center}\begin{tikzpicture}
\matrix(m)[ampersand replacement=\&, matrix of math nodes, row sep=3em, column sep=2.5em, text height=1.5ex, text depth=0.25ex]
{\overline{\hom}(B,\blank)\&\overline{\hom}(A,\blank) \&F\&0\\}; 
\path[->,thick, font=\scriptsize]
(m-1-1) edge node[auto]{}(m-1-2)
(m-1-2) edge node[auto]{}(m-1-3)
(m-1-3) edge node[auto]{}(m-1-4);
\end{tikzpicture}\end{center}   The equivalence is given by  $F\mapsto 	F\circ Q$ where $Q$  is the quotient functor $Q\:\A\to\overline{\A}$.   (2) is dual to (1).  (3) follows by composition of the equivalences in (1) and (2).  The details are left to the reader.  $\qed$\end{pf}

We now state the Serre Localization properties of the restrictions of $\res$ and $\agj$.  The reason these restrict to localizations in the following settings comes from the fact that $w\cong \ev_R\agj$, $w\agl\cong \ev_R$, and $v\cong \ev_R\res$.  The details are left to the reader.

\begin{prop}[Serre Localization Property of $\res$]The restriction of the functor $\res$ to  $$\res\:\fp_0(\Mod(R)^{op},\ab)\to (\underline{\mod}(R)^{op},\ab)$$ is a Serre localization.  That is given any exact functor $$\E\:\fp_0(\Mod(R),\ab)\to \B$$ with $\B$ abelian such that $\ker(\E)\sub \ker(\res)$, then  there is a unique exact functor $\Psi$ making the following diagram commute.\begin{center}\begin{tikzpicture}
\matrix(m)[ampersand replacement=\&, matrix of math nodes, row sep=3em, column sep=3em, text height=1.5ex, text depth=0.25ex]
{\fp_0(\Mod(R),\ab)\&(\underline{\mod}(R)^{op},\ab)\\
\B\&\\};
\path[->,thick, font=\scriptsize]
(m-1-1) edge node[left ]{$ \E$} (m-2-1)
(m-1-2) edge node[auto]{$\exists !\Psi$} (m-2-1)
(m-1-1) edge node[auto]{$\res$}(m-1-2);
\end{tikzpicture}\end{center}\end{prop}

\begin{prop}[Serre Localization Property of $\agj$]The restriction of the Auslander-Gruson-Jensen functor $\agj$ to  $$\agj\:\fp_0(\Mod(R),\ab)\to (\underline{\mod}(R^{op}),\ab)$$ is a Serre localization.  That is for any exact functor $$\E\:\fp_0(\Mod(R),\ab)\to \B$$ with $\B$ abelian and $\ker(\E)\sub \ker(\agj)$, there is a unique exact functor $\Psi$ making the following diagram commute.\begin{center}\begin{tikzpicture}
\matrix(m)[ampersand replacement=\&, matrix of math nodes, row sep=3em, column sep=3em, text height=1.5ex, text depth=0.25ex]
{\fp_0(\Mod(R),\ab)\&(\underline{\mod}(R^{op}),\ab)\\
\B\&\\};
\path[->,thick, font=\scriptsize]
(m-1-1) edge node[left ]{$ \E$} (m-2-1)
(m-1-2) edge node[auto]{$\exists !\Psi$} (m-2-1)
(m-1-1) edge node[auto]{$\agj$}(m-1-2);
\end{tikzpicture}\end{center}\end{prop}

%\begin{pf}This follows from the fact that the functors $F$ and $\agl\agj(F)$ have the same defect and hence $w(F)=0$ if and only if $w(\agl\agj(F))=0$\end{pf}

To state the next result, we need to review the concept of the transpose of a finitely presented module $M$ originally introduced by Auslander in \cite{coh} and later studied by Auslander and Bridger in \cite{stab}.    Suppose that $M\in \mod(R)$.   Then $M$ has presentation

\begin{center}\begin{tikzpicture}
\matrix(m)[ampersand replacement=\&, matrix of math nodes, row sep=3em, column sep=2.5em, text height=1.5ex, text depth=0.25ex]
{R^m\&R^n\&M\&0\\}; 
\path[->,thick, font=\scriptsize]
(m-1-1) edge node[auto]{$A$}(m-1-2)
(m-1-2) edge node[auto]{}(m-1-3)
(m-1-3) edge node[auto]{}(m-1-4);
\end{tikzpicture}\end{center}  where $A$ is an $n$ by $m$ matrix with entries from the ring $R$.  A \dbf{transpose} of the module $M$ is any module $\tr(M)$ obtained by taking the transpose of the matrix $A^T$ and taking its cokernel.  This new module $\tr(M)$ lives in $\mod(R^{op})$ and has a projective presentation
\begin{center}\begin{tikzpicture}
\matrix(m)[ampersand replacement=\&, matrix of math nodes, row sep=3em, column sep=2.5em, text height=1.5ex, text depth=0.25ex]
{R^n\&R^m\&\tr(M)\&0\\}; 
\path[->,thick, font=\scriptsize]
(m-1-1) edge node[auto]{$A^T$}(m-1-2)
(m-1-2) edge node[auto]{}(m-1-3)
(m-1-3) edge node[auto]{}(m-1-4);
\end{tikzpicture}\end{center}  The process of moving from $M$ to a transpose $\tr(M)$ is not functorial.  It depends on the choice of matrix presentation; however, as shown in \cite{stab}, transposes are uniquely determined up to projective equivalence.  That is, if $\tr_1(M)$ and $\tr_2(M)$ are two transposes of $M$, then there exists finitely projectives $P,Q$, such that $$\tr_1(M)\oplus P\cong\tr_2(M)\oplus Q$$  For this reason, if $M\in \mod(R^{op})$ and $F\:\mod(R)\to \ab$ vanishes on projectives, then the assignment $M\mapsto F(\tr(M))$ is actually functorial.  

Recall that a functor $F\:\mod(R)\to \ab$ vanishes on a projective if and only if it factors uniquely through the projectively stable category $\underline{\mod}(R)$ whose objects are the objects in $\mod(R)$ and whose morphisms are morphisms in $\mod(R)$ modulo those factoring through projectives.  Hence there is an equivalence of the functor categories $(\underline{\mod}(R),\ab)$ and the full subcategory of $(\mod(R),\ab)$ consisting of those functors which vanish on projectives.    The functors $F\in (\mod(R),\ab)$ which vanish on projectives are precisely those for which $F(R)=0$.  This means the categories $\ker(\ev_R)$ and $(\underline{\mod}(R),\ab)$ are equivalent.  

%If $F$ vanishes on projectives, then $F(R)=0$ because $R$ is projective.  Assume that $F(R)=0$ and let $P$ be a finitely generated projective in $\mod(R)$.  Since $P$ is a direct summand of $R^n$ for some finite $n$, $F(P)$ is a direct summand of $F(R^n)=F(R)^n=0^n=0$ and hence $F(P)=0$.\end{pf}

In \cite{herzogcontra}, Herzog shows that there is a duality 

\begin{center}\begin{tikzpicture}
\matrix (m) [ampersand replacement= \&,matrix of math nodes, row sep=5em, 
column sep=4em,text height=1.5ex,text depth=0.25ex] 
{(\underline{\mod}(R^{op}),\ab)\&(\underline{\mod}(R)^{op},\ab)\\}; 
\path[->,thick, font=\scriptsize]
(m-1-1)edge[bend right=50] node[below]{$\tr_*$}(m-1-2)
(m-1-2) edge[bend right=50] node[above]{$\tr_*$}(m-1-1);
%(m-1-2) edge[bend left] node[below]{$\agr$}(m-1-1);
\end{tikzpicture}\end{center} given by $\tr_*(F)(M):=F(\tr(M))$.  We can recover this same duality using the functors $\extd$, $\agj$, and $\res$ and the universal property of Serre localization as follows.

\begin{thm}The duality $\extd$ induces a duality $\newdual$ resulting in the following commutative diagram

\begin{center}\begin{tikzpicture}
\matrix(m)[ampersand replacement=\&, matrix of math nodes, row sep=2em, column sep=2em, text height=1.5ex, text depth=0.25ex]
{\fp_0(\Mod(R),\ab)\&\fp_0(\Mod(R)^{op},\ab)\&\fp_0(\Mod(R),\ab)\\
(\underline{\mod}(R^{op}),\ab)\&(\underline{\mod}(R)^{op},\ab)\&(\underline{\mod}(R^{op}),\ab)\\}; 
\path[->,thick, font=\scriptsize]
(m-1-1) edge node[above]{$\extd$}(m-1-2)
(m-1-1) edge node[left]{$\agj$}(m-2-1)
(m-1-2) edge node[above]{$\extd$}(m-1-3)
(m-1-2) edge node[right]{$\res$}(m-2-2)
(m-1-3) edge node[right]{$\agj$}(m-2-3)
(m-2-1) edge node[below]{$\newdual$}(m-2-2)
(m-2-2) edge node[below]{$\newdual$}(m-2-3);
\end{tikzpicture}\end{center} The duality $\newdual$ is isomorphic to the generalized Auslander-Bridger transpose $\tr_*$ studied by Herzog in \cite{herzogcontra}.  Hence the duality $\tr_*$ can be recovered via a universal property.

%\begin{center}\begin{tikzpicture}
%\matrix (m) [ampersand replacement= \&,matrix of math nodes, row sep=5em, 
%column sep=4em,text height=1.5ex,text depth=0.25ex] 
%{\fp_0(\Mod(R),\ab)\&\fp_0(\Mod(R)^{op},\ab)\\}; 
%\path[->,thick, font=\scriptsize]
%(m-1-1)edge[bend right=50] node[below]{$\extd$}(m-1-2)
%(m-1-2) edge[bend right=50] node[above]{$\extd$}(m-1-1);
%%(m-1-2) edge[bend left] node[below]{$\agr$}(m-1-1);
%\end{tikzpicture}\end{center} induces a duality of the Serre quotients
%
%\begin{center}\begin{tikzpicture}
%\matrix (m) [ampersand replacement= \&,matrix of math nodes, row sep=5em, 
%column sep=4em,text height=1.5ex,text depth=0.25ex] 
%{(\underline{\mod}(R^{op}),\ab)\&(\underline{\mod}(R)^{op},\ab)\\}; 
%\path[->,thick, font=\scriptsize]
%(m-1-1)edge[bend right=50] node[below]{$\tr_*$}(m-1-2)
%(m-1-2) edge[bend right=50] node[above]{$\tr_*$}(m-1-1);
%%(m-1-2) edge[bend left] node[below]{$\agr$}(m-1-1);
%\end{tikzpicture}\end{center} 
\end{thm}

\begin{pf}We have the following commutative diagram of functors

\begin{center}\begin{tikzpicture}
\matrix(m)[ampersand replacement=\&, matrix of math nodes, row sep=2em, column sep=2em, text height=1.5ex, text depth=0.25ex]
{\fp_0(\Mod(R),\ab)\&\fp_0(\Mod(R)^{op},\ab)\&\fp_0(\Mod(R),\ab)\\
(\underline{\mod}(R^{op}),\ab)\&(\underline{\mod}(R)^{op},\ab)\&(\underline{\mod}(R^{op}),\ab)\\}; 
\path[->,thick, font=\scriptsize]
(m-1-1) edge node[above]{$\extd$}(m-1-2)
(m-1-1) edge node[left]{$\agj$}(m-2-1)
(m-1-2) edge node[above]{$\extd$}(m-1-3)
(m-1-2) edge node[right]{$\res$}(m-2-2)
(m-1-3) edge node[right]{$\agj$}(m-2-3)
(m-2-1) edge node[below]{$\newdual$}(m-2-2)
(m-2-2) edge node[below]{$\newdual$}(m-2-3);
\end{tikzpicture}\end{center} where $\newdual$ is induced by the Serre localization property of $\agj$ and $\res$ respectively because $\extd$  sends functors arising from pure exact sequences to functors arising from pure exact sequences.  By the uniqueness of $\newdual$, it follows that $\newdual^2\cong 1$ in both cases.  One can easily show that since $\agj((X,\blank))\cong X\otimes\blank$ and $\agj$ is exact, we have $$\agj(\ext^1(M,\blank))\cong \tor_1(M,\blank)$$   Even clearer is that $$\res(\ext^1(\blank, M))\cong \ext^1(\blank, M).$$  Any $F\in(\underline{\mod}(R^{op}),\ab)$ has a projective presentation 
\begin{center}\begin{tikzpicture}
\matrix(m)[ampersand replacement=\&, matrix of math nodes, row sep=3em, column sep=2.5em, text height=1.5ex, text depth=0.25ex]
{\bigoplus_{i\in I}\underline{\hom}(B,\blank)\&\bigoplus_{j\in J}\underline{\hom}(A,\blank)\&F\&0\\}; 
\path[->,thick, font=\scriptsize]
(m-1-1) edge node[auto]{}(m-1-2)
(m-1-2) edge node[auto]{}(m-1-3)
(m-1-3) edge node[auto]{}(m-1-4);
\end{tikzpicture}\end{center} where $A,B$ are finitely presented.  Both $\tr_*$ and $\newdual$ are dualities and hence they are both fully faithful and both send arbitrary sums to products.  Therefore to show that $\tr_*\cong \newdual$ we need only show that they agree on functors of the form $\underline{\hom}(M,\blank)$ where $M$ is finitely presented.   Since $\tr$ is itself a duality on the projectively stable finitely presented module categories, we actually have $$\tr_*(\underline{\hom}(M,\blank))\cong\underline{\hom}(M, \tr(\blank))\cong\underline{\hom}(\blank,\tr(M))$$ Auslander shows in \cite{stab} that  $\underline{\hom}(M,\blank)\cong \tor_1(\tr(M),\blank)$ which we already observed is $\agj(\ext^1(\tr(M),\blank))$.  

\begin{eqnarray*}\newdual\big(\underline{\hom}(M,\blank)\big)&\cong&\newdual\big(\agj(\ext^1(\tr(M),\blank))\big)\\
&\cong&\res\big(\extd\ext^1(\tr(M),\blank)\big)\\
&\cong&\res\bigg(\underline{\hom}\big(\blank, \tr(M)\big)\bigg)\\
&\cong&\underline{\hom}(\blank, \tr(M))
\end{eqnarray*}

This shows that $\newdual$ and $\tr_*$ agree on functors of the form $\underline{\hom}(M,\blank)$ from which it follows that they agree everywhere.   Therefore both $\newdual$ and $\tr_*$ fill the diagram and hence $\newdual\cong \tr_*$ as claimed. $\qed$
\end{pf}
\begin{ex*}[Auslander-Reiten Formulas]
As an example of the Generalised Auslander-Reiten Formulas, we derive the well-known Auslander-Reiten Formulas that appear in the representation theory of Artin algebras, as they are given at \cite{kraus}, from the above theorem. 

Suppose $R$ is a $k$-algebra and let $I$ be an injective cogenerator of $\Mod(k)$, for any commutative ring $k$.

For a functor $F:\A\to\ab$, where $\A$ is a $k$-linear category, we write $D_k(F)=\nat(F\blank,I)$. For a left $R$-module $M$, we also write $D_k(M)$ for the right $R$-module $\hom_k(M,I)$. For functors $F:\A\to\ab$ and $G:A^{op}\to\ab$, there is a natural isomorphism $\nat(F,D_k(G))\cong(G,D_k (F))$ (see \cite[Lemma 2.3]{sam} for a discussion on this isomorphism for modules, which is easily generalised).

One can easily show that $D_k\tor_1(\blank,C)\cong\ext^1(C,D_k\blank)$. This follows easily from the hom-tensor duality, and the fact that $M\mapsto D_kM$ is exact.

In the notation of the above theorem, the Auslander-Reiten formulas become
\begin{align*}
D_k\res\underline{\hom}(\blank,C) &\cong\ext^1(C,D_k\tr\blank)\\
\overline{\hom}(A,D_k\blank)&\cong D_k\newdual\res\ext^{1}(\blank,C)
\end{align*}
We now prove these formulas by using The Generalised Auslander-Reiten formulas and the previous result.

For the first formula,
\begin{align*}
D_k\res\underline{\hom}(\blank,C)&\cong D_k\res\extd\ext^1(C,\blank)
\\&\cong D_k\newdual \agj\ext^1(\blank,C)\\&\cong D_k\newdual\tor_1(\blank,C)\\&\cong D_k\tor_1(\tr\blank,C)\\&\cong\ext^1(C,D_k\tr\blank)
\end{align*}
For the second formula, we need to recall some facts.

For any left $R$-module $M$, $D_kM$ is pure injective. This is a well-known fact (see, e.g. \cite[4.3.29]{psl}). However, it can be seen due to the fact that, in $(\mod(R),\ab)$, $D_k(\blank,M)\cong (D_kM)\otimes_R \blank$, and so, since $(\blank,M)$ is flat in this functor category, its dual is injective. Since $\agl$ agrees with $\agr$ on injectives, it follows that, for any left $R$-module $M$,
\begin{displaymath}
\agl ((D_kM)\otimes\blank)\cong(D_kM,\blank).
\end{displaymath}

From this, we obtain the isomorphisms
\begin{align*}
D_k((\agj\overline{\hom}(A,-))(B))&\cong\nat(\hom(\blank,B),D_k \agj\overline{\hom}(A,\blank))\\
&\cong \nat(\agj\overline{\hom}(A,\blank),D_k\hom(\blank,B))\\
&\cong\nat(\agj\overline{\hom}(A,\blank),(D_kB)\otimes\blank)
\\&\cong \nat(\agl((D_kB)\otimes\blank),\overline{\hom}(A,\blank))
\\&\cong \nat((D_k B,\blank),\overline{\hom}(A,\blank))
\\&\cong\overline{\hom}(A,D_k B)
\end{align*}
That is, 
\begin{displaymath}
D_kD_A\overline{\hom}(A,-)\cong \overline{\hom}(A,D_k\blank).
\end{displaymath}
We can now derive the second Auslander-Reiten formula.
\begin{align*}
\overline{\hom}(A,D_k\blank)&\cong D_k\agj\overline{\hom}(A,-)\\&\cong D_k\agj \extd\ext^1(-,C)
\\&\cong D_k\newdual \res\ext^1(\blank,C)
\end{align*}
\end{ex*}

\begin{rem}It should be pointed out that the isomorphism  $$\nat(\agj\overline{\hom}(A,\blank),(D_kB)\otimes\blank)
\cong \nat(\agl((D_kB)\otimes\blank),\overline{\hom}(A,\blank))$$ used in the example above follows from the adjunction formula of $\agl$ and $\agj$ appearing in \cite{dr1}. \end{rem}

\bibliographystyle{amsplain}
\bibliography{slayer}

\end{document}